\documentclass[letterpaper,journal]{IEEEtran}
\usepackage{amsmath,amsfonts}
\usepackage{array}
\usepackage[caption=false,font=normalsize,labelfont=sf,textfont=sf]{subfig}
\usepackage{textcomp}
\usepackage{stfloats}
\usepackage{url}
\usepackage{verbatim}
\usepackage{graphicx}
\usepackage{cite}
\usepackage{amsmath}             
\usepackage{mathtools}
\usepackage{mathdots}    
\usepackage{amssymb}   
\usepackage{algpseudocode}

\usepackage{graphicx}
\usepackage[table]{xcolor}
\usepackage{multirow}
\usepackage{longtable}
\usepackage{tabularx,booktabs}

\usepackage[normalem]{ulem}

\usepackage{array}

\usepackage[ruled,vlined,linesnumbered]{algorithm2e}

\usepackage{tabularx}
\usepackage{multirow}
\usepackage{booktabs}     
\usepackage[table]{xcolor}  

\begin{document}



\title{EVTradeMatch: A Mobility-Aware Multi-Objective Matching Framework for EV--EV Energy Trading}



\author{
    Md. Mahfujur Rahman,~\IEEEmembership{Graduate Student Member,~IEEE},
    Alistair Barros,
    Raja Jurdak,~\IEEEmembership{Senior Member,~IEEE},
    Darshika Koggalahewa
    
    \thanks{Md Mahfujur Rahman, Alistair Barros, and Darshika Koggalahewa are with the School of Information Systems, Queensland University of Technology, Australia (e-mail: mdmahfujur.rahman@hdr.qut.edu.au; alistair.barros@qut.edu.au; darshika.koggalahewa@qut.edu.au).}
    \thanks{Raja Jurdak is with the School of Computer Science, Queensland University of Technology, Australia (e-mail: r.jurdak@qut.edu.au).}


\thanks{Manuscript received August 11, 2025; revised August 26, 2025.}}




\maketitle




\begin{abstract}

Peer-to-peer energy trading among electric vehicles (EVs) is a promising approach for improving supply-side capacity for accelerated EV adoption and limited charging infrastructure. However, effective EV--EV trading requires coordinated matching between providers and consumers under journey-specific conditions. EV drivers follow distinct journeys and routes with varying traffic densities, arrive at charging locations asynchronously, and have varying energy requirements. Further complexity arises from the dynamic nature of trading partner availability, station capacity, and charging speed across candidate charging locations. 
To address these challenges, we propose EVTradeMatch, a prediction-guided multi-objective optimization framework for mobility-aware EV--EV energy trading. Building on our prior EVNextTrade study for charging-node ranking, we define charging-node suitability as a prediction-derived score that reflects the expected appropriateness of assigning a provider--consumer pair to a candidate charging node, based on mobility, energy, and contextual trading features. This suitability score is used as a guidance signal and as an explicit optimization objective rather than as a hard selection rule.
The EV--EV matching problem is formulated as a multi-objective mixed-integer linear program (MILP) that jointly maximizes matching coverage, transferred energy, and charging-node suitability while minimizing mobility cost under spatial, temporal, and charging-node capacity constraints.
To efficiently explore competing objectives in dynamic, wide-area settings, we develop a tailored non-dominated sorting genetic algorithm II (NSGA-II) within the framework. 
Experimental results show that EVTradeMatch improves transferred energy by 74.2--82.8\% and charging-node suitability by 8.3--84.4\% compared with proximity- and auction-based state-of-the-art methods, while also improving matching coverage by 3.74--25.07 percentage points. The balanced NSGA-II solutions achieve 53.07$\pm$0.76\% matching coverage and transfer 1701.94$\pm$17.03~kWh of energy, with higher mobility cost appearing as an explicit trade-off against methods that prioritize travel minimization. Pareto-front analysis further demonstrates flexible trade-off selection according to operational needs. 
The results demonstrate that combining prediction-derived suitability with multi-objective evolutionary search improves energy throughput and prediction alignment while maintaining competitive matching coverage under the studied journey-aware setting, outcomes not jointly attained by the proximity- or auction-based state-of-the-art methods evaluated here. 

\end{abstract}


\section{Introduction}
\IEEEPARstart{T}{he} rapid growth of electric vehicles (EVs) has intensified pressure on existing charging infrastructures, particularly under dynamic urban mobility patterns and uneven spatial-temporal demand. While public charging stations (CSs) remain the primary mechanism for energy replenishment, their limited availability, fixed deployment, and congestion during peak hours often lead to waiting times, reduced service reliability, and increased range anxiety among EV drivers. In the face of these challenges, new opportunities are emerging to 
address resource limitations in energy supply and distribution in which EV trading is prominent \cite{C62}. In particular, EV–EV energy trading enables peer-to-peer energy exchange between mobile EVs, offering increased flexibility, improved utilization of surplus energy, and enhanced resilience of charging ecosystems \cite{P1, C58}.

Despite its potential, EV–EV energy trading has substantial coordination and optimization challenges \cite{rahman2025survey}. Unlike conventional EV-to-charging-station (EV--CS), EV–EV trading requires matching both consumers and providers, each with independent mobility contexts \cite{li2019intelligent}. Mobile EVs undertake trips for diverse purposes, each with routes having varying traffic densities, and arrive at charging locations at different times, and experience diminishing energy levels during travel \cite{qi2026trustworthy}. In addition, charging locations vary in queuing delay, trading capacity, and rate-constrained energy transfer \cite{huang2019optimal, alinejad2022charge}. Drivers’ willingness to trade may also depend on the area and time of day, making some locations more suitable for trading than others \cite{rahman2025survey}. EV matching, therefore, entails more than simple distance minimization. Together, these mobility, energy, timing, and infrastructure factors make EV--EV matching a constrained multi-objective optimization problem, where improving matching coverage and transferred energy may increase detour distance and waiting time at capacity-limited charging locations.

Existing EV--EV matching studies can be broadly positioned into proximity-based, market-based, and heuristic optimization approaches. Proximity-based methods, such as Gale--Shapley stable marriage \cite{P316} and bichromatic mutual nearest neighbor (BMNN) \cite{P-T-129}, provide efficient matching by emphasizing spatial closeness or pairwise stability. These methods are useful when distance is the dominant matching criterion, yet they are less suitable for synchronous EV--EV trading where feasibility also depends on temporal alignment, energy-transfer capability, charging-location capacity, and journey-specific suitability. Market-based mechanisms, such as double auction matching algorithm (DAMA) \cite{C17}, position EV--EV trading as a bid--ask clearing problem and aim to improve economic efficiency. However, such methods usually abstract away detailed mobility constraints and do not directly optimize the conflicting operational objectives required for provider--consumer coordination in dynamic urban settings. Therefore, the research gap is not simply the absence of matching algorithms, but the lack of a prediction-guided multi-objective matching framework that jointly considers mobility feasibility, temporal coordination, transferred energy, and charging-location suitability.

Shurrab et al. \cite{P316-2} proposed an EV--EV energy trading framework that incorporates multi-parameter optimization, aiming to minimize both cost and charging time while ensuring social welfare and user satisfaction. They used the greedy algorithm and Gale-Shapley matching technique to match providers and consumers based on charging amount and cost. However, the study does not provide details on how their objectives were optimized, and the implementation of mobility (including traffic density) is not provided. It also provides limited research on time and space-constrained optimization, hindering effective coordination in energy trading systems \cite{P-SD}. For instance, energy trading opportunities depend on the spatial proximity of vehicles and their temporal availability, as providers and consumers are required to be co-located to complete a trading transaction. Greedy and rule-based heuristic methods provide scalable solutions by prioritizing immediate gains; however, they lack the ability to explore trade-offs among competing objectives and often lead to suboptimal system-level performance.

In operational EV--EV trading settings, a single aggregated matching solution may not be sufficient because system priorities can change over time. Prior EV charging and e-mobility optimization studies show that charging coordination involves competing objectives across EV users, charging operators, and infrastructure conditions \cite{chung2018charge, hossain2025cloud, song2021multi, cao2023holistic}. During peak charging demand, operators may prioritize maximizing successful provider--consumer matches or transferred energy. During traffic congestion or limited driver availability, reducing mobility cost and waiting time may become more important, as EV routing and charging studies show that travel time, queueing delay, and charging waiting time strongly affect charging decisions \cite{schoenberg2022reducing}. In recommendation-sensitive settings, operators may prefer matches at charging nodes with higher predicted suitability. These changing priorities motivate optimization methods that can generate Pareto-efficient alternatives rather than forcing all objectives into one fixed utility function \cite{deb2002nsga}.

In parallel, recent studies have employed learning-based models to predict charging demand \cite{zafar2025enhancing}, energy availability, or available charging locations \cite{hecht2021predicting}. These prediction tasks improve situational awareness by capturing mobility patterns, energy states, and contextual factors. However, in most existing work, prediction outputs are either used heuristically or are not systematically integrated into downstream constrained matching and multi-objective optimization processes, limiting their effectiveness in coordinated EV--EV energy trading.

In our prior study~\cite{rahman2025evnexttrade}, we developed EVNextTrade, a learning-to-rank framework that ranks candidate WCPs for EV–EV energy trading in journey contexts. Building on this foundation, the present paper addresses a different but complementary problem: network-wide provider--consumer matching for EV--EV energy trading. Specifically, EVTradeMatch takes the prediction-derived suitability scores from EVNextTrade as input and jointly determines which consumer-side EV decision event, provider-side EV decision event, WCP, and time slot should be matched under spatial, temporal, capacity, and one-to-one matching constraints. Therefore, the present study is not merely an extension of charging node recommendation. It transforms recommendation outputs into feasible, multi-objective provider--consumer trading decisions at the network level.



Specifically, we formulate the EV–EV matching problem as a mixed-integer linear program (MILP), because the problem involves discrete matching and trading-location assignment decisions that must satisfy linear constraints on one-to-one matching, spatial reachability, temporal alignment, and station capacity. This formulation provides a rigorous representation of the feasible solution space for provider–consumer matching and trading-location assignment. However, the presence of multiple conflicting objectives and the combinatorial growth of the search space for large-scale dynamic scenarios limit the practicality of exact solution methods alone. We therefore propose a prediction-guided non-dominated sorting
genetic algorithm II (NSGA-II) algorithm, which is well suited to approximating diverse Pareto-efficient solutions without forcing the objectives into a single aggregated utility. The prediction outputs are incorporated as fixed guidance signals to steer the evolutionary search toward more promising provider–consumer matches at feasible trading locations, while ensuring that all final matching decisions remain subject to strict operational constraints.

To ensure fair evaluation, we compare the proposed approach against three baseline methods: random matching, greedy optimization, and PPO (Proximal Policy Optimization, a deep reinforcement learning method) as well as two state-of-the-art methods: BMNN and DAMA-based double auction, all operating on the same feasible matching set. Experimental results demonstrate that prediction-guided multi-objective optimization achieves higher transferred energy and charging-node suitability than conventional approaches, while maintaining competitive matching coverage and making the mobility-cost trade-off explicit.

The main contributions of this paper are summarized as follows:
\begin{itemize}
    \item We provide decision support for EV users in synchronous EV--EV energy trading by identifying feasible provider--consumer matches in journey contexts and assigning them to candidate charging nodes and time slots. The resulting assignments satisfy spatial reachability, temporal compatibility, one-to-one matching, and charging-node capacity constraints, while predicted charging-node suitability is used to prioritize promising trading opportunities.

    
    \item We develop a tailored NSGA-II, a prediction-guided multi-objective optimization algorithm that achieves a balanced exploration of Pareto-optimal trade-offs under strict feasibility requirements.
    
    \item The proposed framework achieves competitive matching performance and balanced trading outcomes for EV--EV matching. When evaluated against heuristic, learning-based, and state-of-the-art approaches, including random matching, greedy optimization, PPO, BMNN, and DAMA, it improves transferred energy and charging-node suitability while maintaining comparable matching coverage, with mobility cost treated as an explicit trade-off rather than as a single-objective minimum.

\end{itemize}

The remainder of this paper is organized as follows. Section 2 reviews related work on synchronous and asynchronous EV–EV trading optimization. Section 3 presents the system model and problem overview. Section 4 introduces the MILP formulation and the NSGA-II algorithm. Section 5 describes the experiment and result analysis, followed by discussion in Section 6. Section 7 summarizes the key observations. Finally, Section 8 concludes the paper and outlines future research directions.

\section{Related Work}

In this section, we briefly review the state-of-the-art studies on optimization relevant to EV-EV trading applications. Optimization is essential for flexible and efficient energy trading. Existing studies on optimization applied to EV applications integrate multiple objectives and constraints to address complex challenges. The objectives include minimizing charging prices, minimizing travel costs (i.e., in terms of travel time and distance), maximizing energy exchanged at the best price for both consumers and providers, maximizing provider-consumer matching, and maximizing EV user social welfare. Social welfare encompasses minimizing costs, reducing wait times, and ensuring fair access to energy resources for the broader EV user community. These objectives are achieved while considering various constraints, such as energy availability, time sensitivity, travel costs, variable traffic conditions, the distance between EVs and CSs, charging speed, and the number of available charging connection ports. 
The existing studies on EV-EV trading optimization can be categorized into two approaches: synchronous trading optimization, which focuses on real-time provider-consumer pair matching, and asynchronous trading optimization, where transactions occur on demand, enabling providers to supply energy to charging storage for later consumption by consumers. As the present study focuses on synchronous EV–EV trading, the following section review existing methodologies proposed in synchronous trading optimization and identify the research gaps that remain open.

\subsubsection{Synchronous trading optimization for provider--consumer pair matching} Existing studies largely focus on optimization of provider-consumer matching while considering different parameters such as distance, balancing supply and demand, cost, charging location etc. to ensure energy consumers and providers participate in EV-EV charging simultaneously. Studies \cite{P316, P-T-129} leveraging the Gale-Shapley stable marriage algorithm and the bichromatic mutual nearest neighbour (BMNN) algorithm have explored distance-based provider-consumer matching involving minimal distances between providers and consumers. However, these approaches have revealed significant limitations. The Gale-Shapley algorithm, designed to create stable pairings by ensuring that no provider and consumer would prefer an alternative match over their current one, may result in suboptimal pairings in complex energy-trading scenarios.
This is because the algorithm prioritizes stability over global optimization, and it has difficulty adapting to unequal sets of providers and consumers. The BMNN algorithm prioritizes distance, overlooking critical factors such as energy availability or time sensitivity, and is less suitable for dynamic environments such as EV mobility, changing energy needs, and variable traffic flow. Both approaches concentrate on local provider-consumer matching within a 5 km range but fail to address mismatches when suitable peers are not found locally. A trip-based probabilistic model for EV-EV energy sharing overcame the limitations of the Gale-Shapley algorithm by considering factors such as energy surplus and detour costs, enabling a more dynamic and context-aware matching process \cite{C11}. All three studies address the coordination aspects of best matching for EV pairing and the distance to charging locations from driving routes.

Many studies have explored optimization techniques for balancing supply and demand and cost optimization in EV charging and energy management \cite{sch2, C19, barros2013exception, sch3, P317}. Wang et al. \cite{sch2} addressed charging overload issues in energy exchange during peak demand hours by optimizing the balance between supply and demand for EVs based on price factors using mixed-integer non-linear programming (MINLP). Similarly, Hosseini et al. \cite{C19} applied the Hungarian algorithm, a robust optimization technique to optimize provider-consumer matching, targeting cost minimization.
However, this approach may inadvertently bias the optimization in favor of consumers, potentially neglecting provider costs and leading to imbalances in energy trading. In addition, as the number of providers and consumers increased, the algorithm was computationally expensive due to its cubic time complexity and additional constraints for matching, making its complexity a bottleneck in large-scale matching problems. To address this scalability issue, He et al. \cite{sch3} extended this optimization process by developing a data-driven, distributionally robust EV balancing method.
By distributing computations across multiple charging nodes in different regions based on EV journeys and energy levels, the proposed approach addressed large-scale optimization problems, mitigating bottlenecks caused by the high number of providers, consumers, and constraints. This method also aimed to minimize the worst-case expected cost while accounting for uncertainties related to passenger mobility demand services, EV supply, and the EV charging state. These studies highlight the importance of sophisticated optimization techniques and pricing factors in enhancing efficiency and fairness in EV energy management, effectively addressing coordination aspects such as charging/discharging control (i.e., supply-demand balance), optimal EV pairing, and large-scale trading support.


The application of bipartite graph algorithms in provider-consumer matching optimization has been explored through various approaches, each aiming to optimize different aspects of the trading process \cite{P317,P-21, P315, li2019intelligent}. Zhang et al. \cite{P317} utilized a weighted bipartite graph to model provider and consumer sets by applying max-weight matching and Gale-Shapley stable marriage algorithms to optimize objectives, such as minimizing costs for individual EV users and maximizing social welfare. Similarly, Zeng et al. \cite{P-21} introduced a hierarchical bipartite graph matching approach to enhance transactive power exchanges among EVs within subdivided distribution power systems. However, their approach overlooked the presence of unmatched EVs, owing to the availability of services. 
Tao et al \cite{P315} also highlighted this issue, emphasizing that existing matching strategies prioritize immediate or short-term efficiency while neglecting long-term benefits and broader implications, such as battery health, cost savings over time, and the overall sustainability of the system. In addition, matching pairs based on bipartite graphs do not account for spatial constraints effectively, leading to impractical solutions where EVs need to travel long distances for charging, which is infeasible in practice \cite{li2019intelligent}. However, both studies have highlighted these limitations without providing a solution. These studies \cite{P317,P-21, P315} considered the coordination aspects of best matching for EV pairing and driver charging station location bias by collectively highlighting the versatility and potential of bipartite graph algorithms in optimizing EV energy trading. However, they had the limitations of handling unmatched EVs, considering long-term benefits, and addressing spatial constraints. Furthermore, most of these studies focus on optimizing a single parameter of EV-EV trading, such as cost, distance, time, social welfare, energy volume, or user preferences.

\begin{figure}[h]
\includegraphics[width=0.5\textwidth]{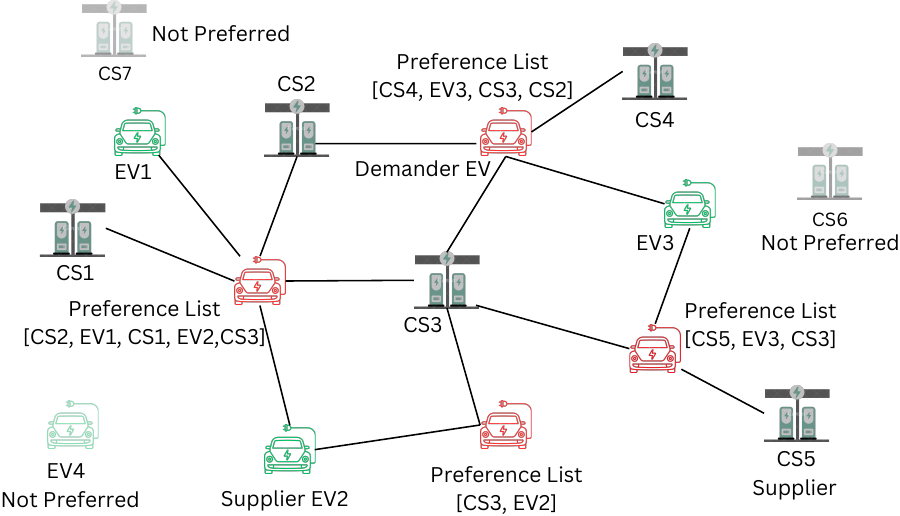}
\caption{Formation of incomplete localized preference lists through local communication \cite{P316}. This figure shows an example scenario where four different EVs request to charge the CS, consumer, and provider EVs in their localized proximity within the same decision time frame. Here, the authors used incomplete terms because consumers will have a partial view of the network.}
\label{fig:options-lit}
\end{figure}

\begin{table*}[!ht]
\footnotesize
\centering
\renewcommand{\arraystretch}{1.4}
\caption{Comparative Study of Existing Research on Synchronous Trading Optimization for Provider--Consumer Pair Matching. Each objective and parameter column indicates whether it is explicitly considered in the corresponding reference. Legends: MonCost: Monetary Cost (charging price, trading revenue, total cost minimization); Dist.: Travel-distance component of mobility cost; Time: Travel-time or waiting-time component of mobility cost; Suit.: Station/Charging-Node Suitability; M$|$T: Mobility$|$Traffic; SoC: State of Charge; Pref.: User Preference/Privacy.}
\label{tab:syn}
\setlength{\tabcolsep}{3pt}
\begin{tabular*}{\textwidth}{@{\extracolsep{\fill}} l p{22mm} *{8}{c} *{5}{c} c @{}}
\hline
\multirow{2}{*}{Ref} & \multirow{2}{=}{\centering Optimization\\Technique} & \multicolumn{8}{c}{Objectives} & \multicolumn{5}{c}{Parameters} & \multirow{2}{*}{M$|$T} \\
\cmidrule(lr){3-10} \cmidrule(lr){11-15}
 & & Match & MonCost & Dist. & Time & Energy & Welfare & Util. & Suit. & Price & Energy & SoC & Loc. & Pref. & \\
\hline
\cite{P316}    & Bipartite + Gale-Shapley         &            &            &            &            &            &            & \checkmark &            &            &            &            &            & \checkmark & X$|$X \\
\hline
\cite{P-T-129} & BMNN                             &            &            & \checkmark &            &            &            &            &            &            &            &            &            & \checkmark & X$|$X \\
\hline
\cite{C19}     & Hungarian algorithm              &            & \checkmark &            &            &            &            &            &            & \checkmark &            &            &            &            & X$|$X \\
\hline
\cite{C10}     & Double-sided auction             &            &            &            &            &            &            & \checkmark &            & \checkmark & \checkmark &            & \checkmark & \checkmark & X$|$X \\
\hline
\cite{P-SD}     & HAPD algorithm                   &            & \checkmark &            &            &            &            & \checkmark &            &            &            &            & \checkmark &            & X$|$X \\
\hline
\cite{P-21}    & Hier. bipartite + PSO            &            &            &            &            & \checkmark &            &            &            &            & \checkmark & \checkmark &            &            & X$|$X \\
\hline
\cite{C17}     & DAMA (double-sided auction)      & \checkmark &            &            & \checkmark & \checkmark &            & \checkmark &            & \checkmark & \checkmark &            &            &            & X$|$X \\
\hline
\cite{sch2}    & Time-coupled MINLP               &            & \checkmark & \checkmark & \checkmark &            &            &            &            & \checkmark & \checkmark &            &            &            & \checkmark$|$X \\
\hline
\cite{P316-2}  & Greedy + Gale-Shapley            &            & \checkmark &            &            & \checkmark & \checkmark &            &            & \checkmark &            & \checkmark &            & \checkmark & \checkmark$|$X \\
\hline
\textbf{Ours}  & \textbf{NSGA-II} & \checkmark &            & \checkmark & \checkmark & \checkmark &            &            & \checkmark &            & \checkmark &            & \checkmark &            & \checkmark$|$\checkmark \\
\hline
\end{tabular*}
\end{table*}

Shurrab et al. \cite{P316-2} have considered multiple parameters for optimization and proposed an EV-EV energy trading framework that aimed at optimizing both cost and charging time, while ensuring social welfare and user satisfaction. They used Greedy algorithm and Gale-Shapley matching technique to match providers and consumers based on charging amount and cost. However, the study does not provide details on how their objectives are optimized, and the implementation of mobility (including traffic density) is not provided. It also limits research on time and space-constrained optimization, hindering effective coordination within energy trading systems \cite{P-SD}. For instance, EVs cannot provide or consume energy arbitrarily in terms of space and time. 

Charging location selection is another optimization parameter for matching providers and consumers, as EVs require a common meeting point to facilitate energy exchange. The core factors in selecting charging locations are trading volume, travel time, charging comfort quality, and charging time \cite{cao2023holistic}. Yassine et al. \cite{C17} proposed integrating edge and fog computing through a double-sided auction mechanism aimed at selecting charging locations while optimizing the energy volume and utility of EV owners from both consumer and provider perspectives. 
In this study, EV users make real-time requests to the system, specifying energy price, energy quantity, and travel time. Based on these inputs, the trading system assigns an appropriate charging location. 
However, their approach assumes that provider and consumer EVs are co-located at any given time, resulting in a static matching. To address this limitation, Ko et al. \cite{P-SD} proposed a software-defined EV-EV provider and consumer matching approach with the support of a mobile charging aggregator. This study aimed to optimize revenue by determining charging locations, costs, and rewards through the implementation of the highly acceptable pair decision (HAPD) algorithm. The study highlighted that matching at fixed charging locations is challenging due to the limited availability of EV charging stations and the high mobility of both provider and consumer EVs. The study proposed the use of dynamic mobile aggregators to facilitate energy transfer from provider to consumer, addressing this issue effectively. These studies address coordination aspects, including charge/supply bias, location bias, charge mode bias, journey bias, and balancing unequal provider-consumer sets for improved energy trading. Table \ref{tab:syn} summarizes the optimization techniques, their objectives, and the parameters for the synchronous EV-EV provider and consumer matching optimization.

Given that most studies focus on synchronous EV-EV trading optimization, the charging duration via this method can be excessively prolonged due to uncertainty regarding EV charging demand or supply. Introducing an asynchronous EV-EV charging service to address these challenges offers a viable solution by decoupling the location and time requirements \cite{C26, cao2023holistic}.

\begin{figure}[h]
\centering
\includegraphics[width=0.5\textwidth]{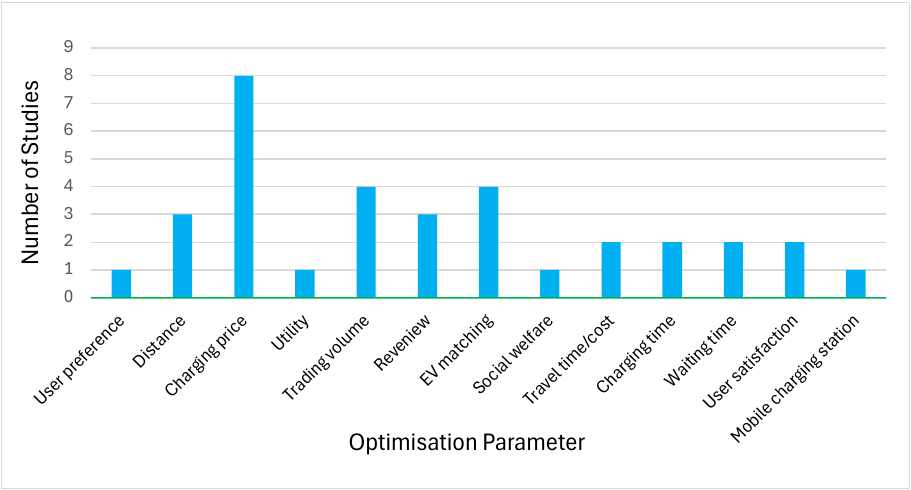}
\caption{This figure illustrates the parameters employed for optimization across various studies.}
\label{fig:options-fig-1}
\end{figure}

\subsubsection{Summary} We found that existing studies often overlook the holistic optimization of mobile EVs in journey contexts within a single study. Studies mainly aim to optimize one or two parameters of EV-EV trading matching, such as charging cost, distance, social welfare, user preference, user satisfaction, required number of mobile charging stations, energy exchange amount, or revenue. Figure \ref{fig:options-fig-1} illustrates the parameters employed for optimization across various studies. 

The comparison in Table \ref{tab:syn} and Figure \ref{fig:options-fig-1} also clarifies the scope of the present study. Monetary cost, price, and social welfare are not included as explicit objectives because this paper focuses on the operational coordination problem of provider--consumer matching and charging-node assignment, rather than market clearing or economic settlement. These economic factors are important in auction-based EV--EV trading, but incorporating them would require additional assumptions on bidding strategies, tariff design, user valuation, and payment mechanisms, which are outside the scope of the proposed mobility-aware matching framework. Similarly, detailed SoC dynamics are not modeled as a standalone optimization objective in this study. Instead, their operational effect is captured through energy demand, provider energy supply, spatial reachability, journey feasibility, and maximum transferable energy constraints. Therefore, the proposed framework targets the spatiotemporal feasibility of synchronous EV--EV trading, while pricing, welfare maximization, and dynamic SoC evolution are left for future extensions.

In addition, several studies state objective functions without clearly specifying the underlying optimization formulation, such as linear programming, mixed-integer programming, or game-theoretic modeling. This limits the reproducibility and interpretability of the proposed approaches, since the decision variables, constraints, and objective trade-offs are often not fully described. Moreover, critical aspects such as temporal feasibility, charging-speed limitations, station capacity, and trade-offs among competing objectives remain underexplored. Based on the reviewed synchronous EV--EV trading literature, no prior study has applied an NSGA-based approach to address the complex trade-offs and constraints involved in provider--consumer matching. This represents an important gap, as NSGA-based optimization can generate a set of Pareto-optimal solutions rather than a single fixed outcome, allowing decision-makers to select the most appropriate matching solution according to the operational context, such as prioritizing higher matching efficiency, lower mobility cost, higher transferred energy, or higher charging-node suitability.

We also found that current optimization approaches often rely on static or pre-specified location assumptions for provider--consumer matching. As a result, they do not sufficiently capture the journey context of mobile EVs, including changing EV locations, changing energy level over time, variable traffic conditions, arrival-time compatibility, and the timeliness of energy transfer. In synchronous EV--EV energy trading, these mobility-dependent factors directly determine whether a provider and consumer can feasibly meet at a charging node, exchange energy within an acceptable time window, and continue their journeys. Therefore, provider--consumer matching should account for mobility-aware and traffic-aware feasibility rather than treating EV location and availability as fixed inputs.

\section{Study Design}
The aim of this study is to develop techniques and algorithms that support efficient EV--EV energy trading by optimally matching provider-side and consumer-side EV decision events at wireless charging pads under dynamic mobility contexts. This requires balancing multiple conflicting objectives, including the number of matched provider--consumer decision-event pairs, transferred energy, mobility cost, and charging-node suitability, while satisfying operational constraints such as spatial reachability, temporal feasibility, and station capacity.

This study is guided by one main research question (RQ): 
\textit{How can provider-side and consumer-side EV decision events be optimally matched at wireless charging pads by balancing multiple conflicting objectives subject to operational constraints, given EV mobility contexts?}

This research question guides the formulation of the EV--EV matching problem, the design of the optimization objectives and constraints, and the evaluation of the proposed framework under mobility-aware trading conditions.

To support this, our contextual assumptions are as follows:

\begin{itemize}
 \item Provider--consumer matching is performed across all candidate EV decision events and charging nodes within the study area, subject to spatial, temporal, and energy feasibility constraints.
 
 
 
 \item Candidate provider--consumer matches are considered feasible only if both EVs can reach the selected charging node and still satisfy their own journey energy requirements after the energy exchange.


  \item EVs have varying energy demand and supply requirements while travelling.
  
 \item Each EV has a source, destination, and possible intermediate journey points.
 \item The distance between EVs and charging pads influences their ability to access energy resources while traveling (i.e., proximity to charging resources).
 \item The decisions made by drivers to determine the best route to their destinations by using navigation applications (e.g., Google Maps or similar services)
 \item Wireless charging pads consider energy sourcing only from provider EVs.
 \item Optimization occurs in one-hour time segments throughout the day.

\end{itemize}

These assumptions will be taken into account to influence the selection of constraints that will be used for optimizing provider-consumer matching. 

\begin{table}[!ht]
\centering
\caption{Notation used in the EV--EV matching formulation.}
\label{tab:notation}
\renewcommand{\arraystretch}{1.15}
\begin{tabularx}{\columnwidth}{p{22mm}X}
\hline
Symbol & Description \\
\hline
$\mathcal{C}$ & Set of consumer EV decision events.\\
$\mathcal{S}$ & Set of provider EV decision events.\\
$\mathcal{K}$ & Set of wireless charging pads (WCPs). \\
$\mathcal{T}$ & Set of discrete time slots in the planning horizon. \\
$K_i \subseteq K$ & Feasible subset of WCPs for consumer EV decision event $i\in\mathcal{C}$, obtained from prediction- and feasibility-based preprocessing. \\

$K_j \subseteq K$ & Feasible subset of WCPs for provider EV decision event $j\in\mathcal{S}$, obtained from prediction- and feasibility-based preprocessing. \\
$d_{ik}$ & Travel distance from consumer $i$ to WCP $k$. \\
$d_{jk}$ & Travel distance from provider $j$ to WCP $k$. \\
$a_{ik}$ & Arrival time of consumer $i$ at WCP $k$. \\
$a_{jk}$ & Arrival time of provider $j$ at WCP $k$. \\
$t_i^{\mathrm{dep}}$ & Departure time of consumer $i$. \\
$t_j^{\mathrm{dep}}$ & Departure time of provider $j$. \\
$W_i$ & Maximum waiting tolerance of consumer $i$. \\
$W_j$ & Maximum waiting tolerance of provider $j$. \\
$\mathrm{Cap}_{k}^{t}$ & Capacity of WCP $k$ in slot $t$, measured as the maximum number of simultaneous EV--EV charging sessions. \\
$E^t_i$ & Energy demand associated with consumer decision event $i$. \\
$E^t_j$ & Energy supply associated with provider decision event $j$. \\
$p_k$ & Charging speed of WCP $k$. \\
$\Delta$ & Duration of each time slot in hours. \\

$s^c_{ik}$ & Raw prediction score obtained from the EVNextTrade recommendation model for consumer $i$ at WCP $k$. \\

$s^p_{jk}$ & Raw prediction score obtained from the EVNextTrade recommendation model for provider $j$ at WCP $k$.\\

$q^c_{ik}$ & Prediction-derived suitability probability of WCP $k$ for consumer $i$, obtained as $q^c_{ik}=\sigma(s^c_{ik})$.
\\
$q^p_{jk}$ & Prediction-derived suitability probability of WCP $k$ for provider $j$, obtained as $q^p_{jk}=\sigma(s^p_{jk})$.
\\
$q_{ijk}$ &  Pairwise suitability score for matching consumer $i$ with provider $j$ at WCP $k$, computed as $q_{ijk}=(q^c_{ik}+q^p_{jk})/2$.\\


$\mathrm{Dist}^{c}_{\max}$ & Maximum allowable detour distance for consumers. \\
$\mathrm{Dist}^{s}_{\max}$ & Maximum allowable detour distance for providers. \\
$A^c_{ik}$ & Binary indicator equal to 1 if consumer $i$ can feasibly reach WCP $k$, and 0 otherwise. \\
$A^s_{jk}$ & Binary indicator equal to 1 if provider $j$ can feasibly reach WCP $k$, and 0 otherwise. \\
$B^c_{ik,t}$ & Binary indicator equal to 1 if consumer $i$ arrives at WCP $k$ during slot $t$, and 0 otherwise. \\
$B^s_{jk,t}$ & Binary indicator equal to 1 if provider $j$ arrives at WCP $k$ during slot $t$, and 0 otherwise. \\
$A_{ijk,t}$ & Binary temporal-feasibility indicator for consumer $i$, provider $j$, WCP $k$, and slot $t$. \\
$X_{ijk}^{t}$ & Binary decision variable equal to 1 if consumer decision event $i$ is matched with provider decision event $j$ at WCP $k$ in slot $t$, and 0 otherwise. \\
$\bar{E}_{ijk}$ & Maximum transferable energy for a feasible match between consumer $i$ and provider $j$ at WCP $k$ within one slot. \\
\hline
\end{tabularx}
\end{table}

\section{Research Methodology}
\label{sec:methodology}


\begin{figure*}[h]
\centering
\includegraphics[width=\textwidth]{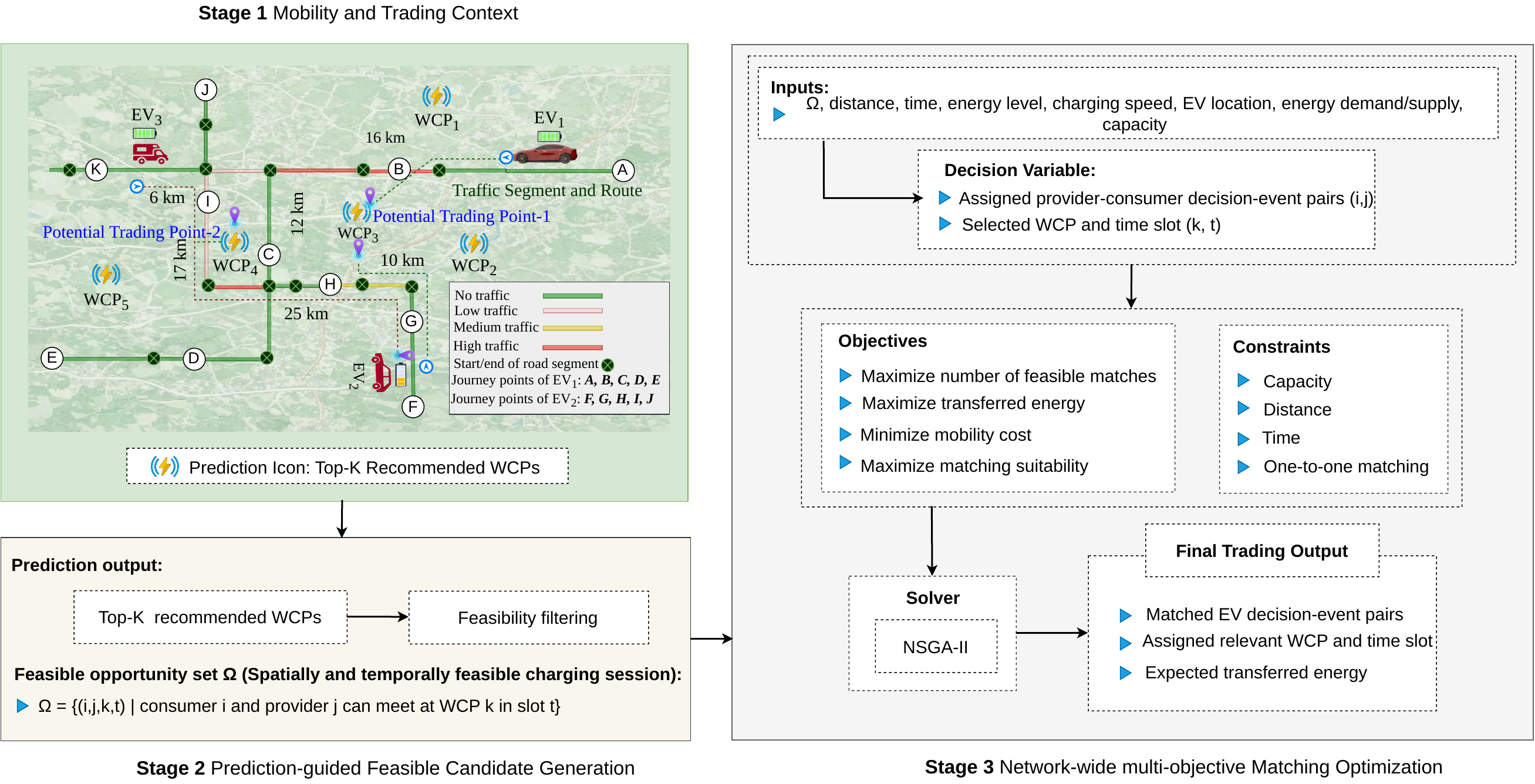}
\caption{Prediction-guided EV–EV energy trading optimization. Prediction outputs identify spatially and temporally feasible charging opportunities, forming the feasible opportunity set $\Omega$. The proposed NSGA-II algorithm selects EV–EV charging sessions from $\Omega$ by optimizing matching coverage, transferred energy, mobility cost, and matching suitability under distance, time, capacity, and one-to-one matching constraints. Capacity is enforced per wireless charging pad and time slot (k,t), while each consumer or provider EV decision event participates in at most one selected charging session} over the optimization horizon.
\label{fig:options-fig-2}
\end{figure*}

\subsection{Overview}

We consider a mobility-aware EV--EV energy trading setting where consumer-side EV decision events (demand) are matched with provider-side EV decision events (supply) at wireless charging pads (WCPs). In this study, an EV decision event denotes a role-specific consumer or provider trading request generated by a physical EV under a particular journey, energy, and candidate-WCP context. The decision-making problem is to construct feasible provider--consumer matches under spatiotemporal constraints induced by mobility and bounded waiting tolerance, as well as infrastructure limitation such as station capacity. In our overall research framework, charging-node suitability scores are obtained from EVNextTrade, a learning-to-rank-based recommendation model developed in our prior study~\cite{rahman2025evnexttrade}. EVNextTrade ranks candidate charging nodes for each EV decision event and produces prediction-derived relevance scores that reflect the relative suitability of each node for EV--EV energy trading. In this paper, these scores are used as fixed inputs to the provider--consumer matching optimization. Because matching suitability depends on both consumer-side and provider-side conditions, we combine the corresponding prediction-derived scores into a pairwise suitability score for each feasible provider--consumer--charging-node tuple. The rank or predicted suitability of candidate charging nodes is therefore incorporated into the matching process as a soft guidance signal and as an explicit suitability objective, rather than as a hard rule that directly selects the top-ranked node. Final matching decisions remain subject to the optimization constraints. This prediction-to-optimization linkage forms the basis of the proposed prediction-guided framework.

\subsection{Sets and Parameters}
\label{sec:notation}

This section introduces the notation used in the proposed EV--EV matching model. The formulation considers consumer EV decision events, provider EV decision events, wireless charging pads (WCPs), and discrete time slots, together with mobility, energy, and infrastructure-related parameters. Table~\ref{tab:notation} summarizes the main symbols used throughout the formulation.

\noindent\textbf{EV decision-event indexing.}
Throughout the formulation and algorithms, the consumer index $i\in\mathcal{C}$ and provider index $j\in\mathcal{S}$ refer to EV decision events. 
An EV decision event denotes a role-specific consumer or provider trading request generated by a physical EV at a particular time, location, energy level, and candidate-WCP context. Therefore, all quantities indexed by $i$ and $j$, including distance, arrival time, energy demand or supply, feasibility indicators, suitability scores, and matching decisions, are defined at the EV decision-event level. 

\noindent\emph{Sets.}
The following sets represent the EV decision events, charging infrastructure, and planning horizon:
\begin{itemize}
    \item $\mathcal{C}$: set of consumer EV decision events.
    \item $\mathcal{S}$: set of provider EV decision events.
    \item $\mathcal{K}$: set of wireless charging pads (WCPs).
    \item $\mathcal{T}$: set of discrete time slots.
    \item $\mathcal{K}_i \subseteq \mathcal{K}$: feasible WCPs for consumer EV decision event $i \in \mathcal{C}$.
    \item $\mathcal{K}_j \subseteq \mathcal{K}$: feasible WCPs for provider EV decision event $j \in \mathcal{S}$.
\end{itemize}

\noindent\emph{Parameters.}
The model uses the following mobility-, energy-, and infrastructure-related parameters:
\begin{itemize}
    \item $d_{ik}$: travel distance from consumer $i$ to WCP $k$.
    \item $d_{jk}$: travel distance from provider $j$ to WCP $k$.
    \item $a_{ik}$: arrival time of consumer $i$ at WCP $k$.
    \item $a_{jk}$: arrival time of provider $j$ at WCP $k$.
    \item $t_i^{\mathrm{dep}}$: departure time of consumer $i$.
    \item $t_j^{\mathrm{dep}}$: departure time of provider $j$.
    \item $W_i$: waiting tolerance of consumer $i$.
    \item $W_j$: waiting tolerance of provider $j$.
    \item $\mathrm{Cap}_{k}^{t}$: capacity of WCP $k$ during slot $t$.
    \item $E^t_i$: energy demand associated with consumer decision event $i$.
    \item $E^t_j$: energy supply associated with provider decision event $j$.
    \item $p_k$: charging speed of WCP $k$.
    \item $\Delta$: slot duration in hours.
    \item $\mathrm{Dist}^{c}_{\max}$: maximum allowable detour distance for consumers.
    \item $\mathrm{Dist}^{s}_{\max}$: maximum allowable detour distance for providers.
    \item $s^c_{ik}$: raw prediction score obtained from the EVNextTrade recommendation model for consumer $i$ at WCP $k$.
\item $s^p_{jk}$: raw prediction score obtained from the EVNextTrade recommendation model for provider $j$ at WCP $k$.
\item $q^c_{ik}$: prediction-derived suitability probability of WCP $k$ for consumer $i$, obtained as $q^c_{ik}=\sigma(s^c_{ik})$.
\item $q^p_{jk}$: prediction-derived suitability probability of WCP $k$ for provider $j$, obtained as $q^p_{jk}=\sigma(s^p_{jk})$.
\item $q_{ijk}$: pairwise suitability score for matching consumer $i$ with provider $j$ at WCP $k$, computed as $q_{ijk}=(q^c_{ik}+q^p_{jk})/2$.
    \item $\sigma_{\epsilon}$: tie-breaking noise scale in the prediction-guided scoring function, used to break ties between candidates with nearly identical utilities and thereby produce diverse individuals in the initial population.
\end{itemize}


\noindent\emph{Feasibility indicators.}
To encode spatial and temporal reachability within the MILP, we precompute the following binary indicators:
\begin{itemize}
    \item $A^c_{ik}\in\{0,1\}$: equals 1 if consumer $i$ can feasibly reach WCP $k$, and 0 otherwise.
    \item $A^s_{jk}\in\{0,1\}$: equals 1 if provider $j$ can feasibly reach WCP $k$, and 0 otherwise.
    \item $B^c_{ik,t}\in\{0,1\}$: equals 1 if consumer $i$ arrives at WCP $k$ in slot $t$, and 0 otherwise.
    \item $B^s_{jk,t}\in\{0,1\}$: equals 1 if provider $j$ arrives at WCP $k$ in slot $t$, and 0 otherwise.
    \item $A_{ijk,t}\in\{0,1\}$: equals 1 if consumer $i$ and provider $j$ can be temporally matched at WCP $k$ in slot $t$, and 0 otherwise.
\end{itemize}

\subsection{Decision Variable}

We define the binary decision variable $X_{ijk}^{t} \in \{0,1\}$ to indicate whether consumer EV decision event $i$ is matched with provider EV decision event $j$ at wireless charging pad $k$ during time slot $t$:
\[
X_{ijk}^{t} =
\begin{cases}
1, & \text{if } i \text{ is matched with } j \text{ at } k \text{ in } t, \\
0, & \text{otherwise.}
\end{cases}
\]

\subsection{Constraints}
\label{sec:cons}

\label{sec:cons}

\noindent (i) One-to-one matching for consumers:  
Each consumer EV decision event can participate in at most one EV--EV charging session within the planning horizon. This prevents a consumer-side decision event from being assigned to multiple providers, charging pads, or time slots.
\begin{equation}
\sum_{t\in\mathcal{T}} \sum_{k\in\mathcal{K}} \sum_{j\in\mathcal{S}} X_{ijk}^{t} \le 1,
\quad \forall i\in\mathcal{C}.
\label{eq:cons_consumer_unique}
\end{equation}

\noindent (ii) One-to-one matching for providers:  
Similarly, each provider EV decision event can be matched to at most one consumer EV decision event within the planning horizon. This ensures that a provider-side decision event cannot be assigned to multiple charging sessions at different locations or times.
\begin{equation}
\sum_{t\in\mathcal{T}} \sum_{k\in\mathcal{K}} \sum_{i\in\mathcal{C}} X_{ijk}^{t} \le 1,
\quad \forall j\in\mathcal{S}.
\label{eq:cons_provider_unique}
\end{equation}

\noindent (iii) WCP capacity per time slot:  
Each wireless charging pad has a limited capacity in each time slot, representing the maximum number of simultaneous EV--EV charging sessions it can support. This constraint ensures that the number of matched provider--consumer pairs assigned to charging pad $k$ in slot $t$ does not exceed its available capacity.
\begin{equation}
\sum_{i\in\mathcal{C}} \sum_{j\in\mathcal{S}} X_{ijk}^{t} \le \mathrm{Cap}_{k}^{t},
\quad \forall k\in\mathcal{K},\; \forall t\in\mathcal{T}.
\label{eq:cons_capacity}
\end{equation}

\noindent (iv) Spatial feasibility:
To ensure that the EVs associated with a consumer decision event and a provider decision event can feasibly reach a selected WCP, their travel distances to that WCP must not exceed their respective maximum allowable detour distances. Let $d_{ik}$ and $d_{jk}$ denote the travel distances of the EVs associated with consumer decision event $i$ and provider decision event $j$ to WCP $k$, respectively. Let $\mathrm{Dist}^{c}_{\max}$ and $\mathrm{Dist}^{s}_{\max}$ denote the maximum allowable detour distances for consumers and providers.

We define the spatial feasibility indicators as
\begin{equation}
A^{c}_{ik} =
\begin{cases}
1, & \text{if } d_{ik} \le \mathrm{Dist}^{c}_{\max}, \\
0, & \text{otherwise,}
\end{cases}
\label{eq:spatial_indicator_consumer}
\end{equation}
and
\begin{equation}
A^{s}_{jk} =
\begin{cases}
1, & \text{if } d_{jk} \le \mathrm{Dist}^{s}_{\max}, \\
0, & \text{otherwise.}
\end{cases}
\label{eq:spatial_indicator_provider}
\end{equation}

The resulting feasibility conditions are enforced through
\begin{equation}
X_{ijk}^{t} \le A^{c}_{ik},
\quad \forall i\in\mathcal{C},\; \forall j\in\mathcal{S},\; \forall k\in\mathcal{K},\; \forall t\in\mathcal{T},
\label{eq:cons_spatial_c}
\end{equation}
and
\begin{equation}
X_{ijk}^{t} \le A^{s}_{jk},
\quad \forall i\in\mathcal{C},\; \forall j\in\mathcal{S},\; \forall k\in\mathcal{K},\; \forall t\in\mathcal{T}.
\label{eq:cons_spatial_s}
\end{equation}


\noindent (v) Temporal feasibility:
To ensure that consumer decision event $i$ and provider decision event $j$ can be paired for an EV--EV charging session at wireless charging pad $k$ within a discrete time slot $t$, their arrival times must be sufficiently close under their waiting tolerances and must fall within the same charging slot.

Let $t_i^{\mathrm{dep}}$ and $t_j^{\mathrm{dep}}$ denote the departure times of the EVs associated with consumer decision event $i$ and provider decision event $j$, respectively. Let $a_{ik}$ and $a_{jk}$ denote their arrival times at WCP $k$, computed as the sum of departure time and route travel time under time-dependent traffic conditions:
\[
a_{ik} = t_i^{\mathrm{dep}} + \sum_{e \in \pi_{ik}} \frac{\ell_e}{v_e(t_i^{\mathrm{dep}})},
\qquad
a_{jk} = t_j^{\mathrm{dep}} + \sum_{e \in \pi_{jk}} \frac{\ell_e}{v_e(t_j^{\mathrm{dep}})},
\]
where $\pi_{ik}$ and $\pi_{jk}$ denote the sets of road segments on the routes from EV $i$ and EV $j$ to WCP $k$, respectively, $\ell_e$ is the length of segment $e$, and $v_e(\cdot)$ is the time-dependent traffic speed on segment $e$. Traffic-dependent travel times are computed using temporally discretized traffic conditions evaluated at the EV departure time, consistent with the granularity of the available traffic data and the charging-slot formulation.

Let $W_i$ and $W_j$ denote the maximum allowable waiting tolerances of EV $i$ and EV $j$. A necessary condition for temporal overlap is
\begin{equation}
|a_{ik} - a_{jk}| \le \min(W_i, W_j).
\label{eq:time_overlap_condition}
\end{equation}



Let $\mathcal{T} = \{0, 1, \ldots, T-1\}$ denote the set of discrete time slots, where each slot $t \in \mathcal{T}$ corresponds to the half-open interval $[t\Delta,(t+1)\Delta)$. This convention ensures that each arrival time belongs to exactly one slot, with no overlap or gaps at slot boundaries. We define the slot-membership indicators as

\begin{equation}
\begin{aligned}
B^c_{ik,t} &=
\begin{cases}
1, & \text{if } a_{ik} \in [t\Delta,(t+1)\Delta), \\
0, & \text{otherwise,}
\end{cases} \\
B^s_{jk,t} &=
\begin{cases}
1, & \text{if } a_{jk} \in [t\Delta,(t+1)\Delta), \\
0, & \text{otherwise.}
\end{cases}
\end{aligned}
\label{eq:slot_membership}
\end{equation}

Using these indicators, the slot-indexed temporal-feasibility indicator is defined as
\begin{equation}
A_{ijk,t} =
\begin{cases}
1, &
\begin{aligned}
& |a_{ik} - a_{jk}| \le \min(W_i, W_j) \\
& \wedge\; B^c_{ik,t}=1 \\
& \wedge\; B^s_{jk,t}=1,
\end{aligned} \\
0, & \text{otherwise.}
\end{cases}
\label{eq:time_indicator_slot}
\end{equation}

The temporal-feasibility constraint is then enforced as
\begin{equation}
X_{ijk}^{t} \le A_{ijk,t},
\quad \forall i\in\mathcal{C},\; \forall j\in\mathcal{S},\; \forall k\in\mathcal{K},\; \forall t\in\mathcal{T}.
\label{eq:cons_time_gate}
\end{equation}



\noindent (vi) Binary domain:
\begin{equation}
X_{ijk}^{t} \in \{0,1\},
\quad \forall i\in\mathcal{C},\; \forall j\in\mathcal{S},\; \forall k\in\mathcal{K},\; \forall t\in\mathcal{T}.
\label{eq:cons_binary}
\end{equation}

\subsection{Model Assumptions}

The optimization is performed over a fixed planning horizon at the EV decision-event level. Within this horizon, each consumer or provider EV decision event is allowed to participate in at most one selected EV--EV trading session. Although a physical EV may generate multiple EV decision events over the planning horizon, these events are treated as distinct role-specific trading requests because they may differ in time, location, energy state, and candidate-WCP context. Time is discretized into slots to capture session-level temporal feasibility and infrastructure usage, while exact scheduling within a selected slot is handled as a deterministic post-processing step.
The current formulation incorporates energy and temporal feasibility at the EV decision-event level. Each selected session must satisfy event-level energy demand or supply, arrival-time compatibility, spatial reachability, and WCP--slot capacity constraints. 

\subsection{Objective Functions}
\label{sec:milp}

We formulate provider--consumer matching as a multi-objective mixed-integer linear program that simultaneously maximizes matching coverage, measured by the number of selected EV--EV trading sessions, transferred energy, and charging-node suitability while minimizing mobility cost.

\noindent (i) Maximize matching coverage, measured as the number of selected EV--EV trading sessions:
\begin{equation}
\max f_1 = \sum_{t\in\mathcal{T}} \sum_{k\in\mathcal{K}} \sum_{i\in\mathcal{C}} \sum_{j\in\mathcal{S}} X_{ijk}^{t}
\label{eq:obj_matches}
\end{equation}

\noindent (ii) Maximize transferred energy:
The maximum transferable energy for a match at WCP $k$ within one slot is limited by demand, supply, and charging speed:
\begin{equation}
\bar{E}_{ijk} = \min\left(E^t_i,\; E^t_j,\; p_k \cdot \Delta\right),
\label{eq:transfer_cap}
\end{equation}
which is treated as a constant during optimization. The second objective is therefore
\begin{equation}
\max f_2 = \sum_{t\in\mathcal{T}} \sum_{k\in\mathcal{K}} \sum_{i\in\mathcal{C}} \sum_{j\in\mathcal{S}} \bar{E}_{ijk}\, X_{ijk}^{t}
\label{eq:obj_energy}
\end{equation}


\noindent (iii) Minimize mobility cost: The mobility cost is defined as a weighted combination of normalized total travel distance and normalized waiting-time mismatch, where $\alpha$ and $\beta$ are weighting coefficients for the distance and waiting-time components, respectively, $D_{\max}$ denotes the maximum allowed total travel distance, and $W_{\max}$ denotes the maximum allowed waiting tolerance.
\begin{equation}
\min f_3 =
\sum_{t\in\mathcal{T}} \sum_{k\in\mathcal{K}} \sum_{i\in\mathcal{C}} \sum_{j\in\mathcal{S}}
\left[
\alpha \frac{d_{ik}+d_{jk}}{D_{\max}}
+
\beta \frac{|a_{ik}-a_{jk}|}{W_{\max}}
\right] X_{ijk}^{t}
\label{eq:obj_mobility}
\end{equation}

\textbf{Definition 1 (Suitability Score).}
For each candidate WCP $k \in \mathcal{K}$, consumer decision event $i \in \mathcal{C}$ and provider decision event $j \in \mathcal{S}$ are associated with raw prediction scores $s_{ik}^{c}, s_{jk}^{p} \in \mathbb{R}$, obtained from the EVNextTrade recommendation model~\cite{rahman2025evnexttrade}. These scores estimate the suitability of consumer $i$ and provider $j$ for participating in energy trading at WCP $k$. The raw scores are transformed into bounded suitability probabilities using the logistic sigmoid function $\sigma(\cdot)$:
\begin{equation}
q_{ik}^{c} = \sigma(s_{ik}^{c}),
\qquad
q_{jk}^{p} = \sigma(s_{jk}^{p}),
\qquad
\sigma(x)=\frac{1}{1+e^{-x}}.
\label{eq:suitability_components}
\end{equation}
The pairwise suitability score for matching consumer $i$ with provider $j$ at WCP $k$ is then defined as:
\begin{equation}
q_{ijk}
=
\frac{q_{ik}^{c}+q_{jk}^{p}}{2},
\qquad
q_{ijk}\in(0,1).
\label{eq:suitability_pair}
\end{equation}
Thus, $q_{ijk}$ represents the combined predicted suitability of the consumer--provider pair $(i,j)$ at WCP $k$.

The arithmetic mean is used to assign symmetric weight to the consumer-side and provider-side suitability scores. This reflects the assumption that a trading session is desirable only when the selected WCP is suitable for both participants: the consumer receiving energy and the provider supplying energy. Thus, $q_{ik}^{c}$ and $q_{jk}^{p}$ each receive a weight of $0.5$ in the pairwise suitability score, avoiding a design bias toward either participant. 



\noindent (iv) Maximize predicted charging-node suitability:
To favor matches that are more suitable for the selected charging node, the fourth objective maximizes the aggregated suitability score of all selected consumer--provider--station--time assignments:
\begin{equation}
\max f_4
=
\sum_{t\in\mathcal{T}}
\sum_{k\in\mathcal{K}}
\sum_{i\in\mathcal{C}}
\sum_{j\in\mathcal{S}}
q_{ijk} X_{ijk}^{t}
\label{eq:obj_suitability}
\end{equation}


\noindent Subject to:
\begin{itemize}
    \item consumer and provider decision-event uniqueness constraints in Eqs.~\eqref{eq:cons_consumer_unique}--\eqref{eq:cons_provider_unique},
    \item charging-pad capacity constraints in Eq.~\eqref{eq:cons_capacity},
    \item spatial feasibility constraints in Eqs.~\eqref{eq:cons_spatial_c}--\eqref{eq:cons_spatial_s},
    \item temporal feasibility constraints in Eq.~\eqref{eq:cons_time_gate}, and
    \item binary decision-variable constraints in Eq.~\eqref{eq:cons_binary}.
\end{itemize}

\subsection{Matching Mechanisms}
\label{sec:pg-nsga2}

This section presents a matching procedure to construct near-optimal provider--consumer decision-event matches at wireless charging pads (WCPs) under mobility-driven spatiotemporal feasibility, bounded waiting tolerance, and infrastructure capacity limits. The search space is combinatorial because each feasible solution must jointly decide (i) which consumer--provider decision-event pairs are formed, (ii) at which WCP they meet, and (iii) in which time slot the EV--EV charging session occurs, while enforcing one-to-one matching and WCP--slot capacity constraints.

In our setting, prediction-derived suitability scores for candidate WCPs are produced by an upstream prediction model and treated as fixed inputs to the matching optimization. Here, \textit{upstream} refers to the prediction model from our prior study~\cite{rahman2025evnexttrade}, which generates candidate-WCP suitability scores prior to the optimization stage. These prediction outputs are combined with spatial-reachability and temporal-alignment checks to restrict the optimization to a feasible set of candidate charging sessions through binary feasibility indicators. In addition, the predicted suitability score for each candidate session is directly incorporated as an explicit objective in the multi-objective optimization.

The matching optimization simultaneously optimizes four objectives: (i) maximizing matching coverage, (ii) maximizing transferred energy, (iii) minimizing mobility cost, and (iv) maximizing predicted suitability of selected charging sessions, as defined in Eqs.~\eqref{eq:obj_matches}--\eqref{eq:obj_suitability}, subject to feasibility and capacity constraints in Eqs.~\eqref{eq:cons_consumer_unique}--\eqref{eq:cons_binary}. Suitability scores are therefore used directly in the optimization, ensuring that the selected matches are not only feasible and efficient but also aligned with the most promising charging nodes identified by the upstream prediction model.

To optimize these competing objectives, we employ NSGA-II with domain-specific enhancements. Specifically, both population initialization and mutation are \emph{prediction-guided}: when constructing or modifying solutions, candidate matches are preferentially sampled from the sets of WCPs that are feasible for each consumer EV decision event ($\mathcal{K}_i$) and provider EV decision event ($\mathcal{K}_j$). Among these feasible candidates, the selection is further guided toward WCPs with higher predicted suitability scores ($q_{ijk}$), meaning that more suitable charging pads are more likely to be chosen. After each construction or mutation step, a repair procedure is applied to ensure that all solutions satisfy one-to-one matching and WCP--slot capacity constraints. The repair procedure checks each candidate solution and modifies it as needed to enforce feasibility: if a consumer or provider EV decision event is assigned to more than one session, or if the number of sessions at any WCP in a time slot exceeds its capacity, lower-scoring conflicting assignments are removed until all one-to-one and capacity constraints are satisfied. The complete procedure is detailed in Algorithms~\ref{alg:pg-nsga2-main}--\ref{alg:evaluate-fitness}.

The proposed MILP formulation provides a rigorous mathematical representation of the EV--EV matching problem, including one-to-one matching, spatial reachability, temporal feasibility, WCP capacity, transferable energy, mobility cost, and prediction-derived suitability. The formulation is used to define the feasible decision space, objective structure, and constraint set. Exact multi-objective MILP solving is not used in the main experiments because the number of feasible provider--consumer--WCP--time-slot tuples grows combinatorially with the number of EV decision events and candidate charging nodes, making exact optimization computationally demanding for realistic problem sizes. Therefore, we use NSGA-II as a scalable multi-objective search method to approximate the Pareto-efficient solution set under the same objectives and constraints. NSGA-II is well suited to this setting because it can handle constrained, non-convex, and discontinuous combinatorial search spaces, maintain solution diversity, and incorporate domain-specific mechanisms such as prediction-guided initialization and feasibility repair.

\subsubsection{\textbf{NSGA-II Approach for Mobility-Aware EV--EV Matching}}
\label{sec:pg-nsga2-overview}

Algorithm~\ref{alg:pg-nsga2-main} presents the overall workflow of the proposed NSGA-II algorithm for mobility-aware EV--EV
matching. The algorithm systematically explores the feasible matching space to identify high-quality trade-off solutions that maximize matching coverage and energy transfer while minimizing mobility-induced cost, subject to one-to-one, capacity, and spatiotemporal feasibility constraints.

\paragraph{Input and Output}
Algorithm~\ref{alg:pg-nsga2-main} accepts as input the number of generations $G$, population size $N$, mutation probability $p_m$, and the precomputed feasible session universe $\Omega$. It also receives the sets of consumers $\mathcal{C}$, providers $\mathcal{S}$, wireless charging pads $\mathcal{K}$, and time slots $\mathcal{T}$; feasibility indicators $A^c_{ik}$, $A^s_{jk}$, and $A_{ijk,t}$; WCP capacities $\mathrm{Cap}_k^t$; and the optimization parameters $q_{ijk}$, $\bar{E}_{ijk}$, $d_{ik}$, $d_{jk}$, $a_{ik}$, $a_{jk}$, $\alpha$, and $\beta$. The repair parameters $\lambda$ and $\mu$ are used by the feasibility-repair procedure. The output is the first non-dominated front $\mathcal{F}^{\star}$, where each solution represents a feasible EV--EV matching configuration consisting of selected sessions $(i,j,k,t)$.

\paragraph{Feasible Session Universe}
The algorithm uses the precomputed feasible session universe $\Omega$ in line~1. Each tuple $(i,j,k,t)\in\Omega$ represents a realizable EV--EV charging session between consumer decision event $i$ and provider decision event $j$ at WCP $k$ during time slot $t$. Since $\Omega$ is constructed before the NSGA-II search using the spatial and temporal feasibility indicators, the evolutionary search is restricted to candidate sessions that satisfy the basic feasibility requirements.


\paragraph{Population Initialization} The initial population is generated in two steps. First, line~2 constructs an intermediate population $\mathcal{P}_{0}^{init}$ of $N$ individuals using \textbf{SelectSession} from Algorithm~\ref{alg:select-session}. Each individual is formed by selecting feasible provider--WCP--time assignments for consumer EV decision events. Then, line~3 applies \textbf{RepairFeasibility} from Algorithm~\ref{alg:repair-feasibility} to each individual in $\mathcal{P}_{0}^{init}$ and stores the repaired population as $\mathcal{P}_{0}$. This ensures that the initial population contains candidate matching solutions that respect the feasible session universe and are corrected for consumer, provider, and WCP-capacity conflicts.

\begin{algorithm*}[t]
\small
\caption{NSGA-II for Mobility-Aware EV--EV Matching}
\label{alg:pg-nsga2-main}
\DontPrintSemicolon
\LinesNumbered

\KwIn{$G$ (number of generations); $N$ (population size); mutation probability $p_m$; precomputed feasible session universe $\Omega$; Sets $\mathcal{C},\mathcal{S},\mathcal{K},\mathcal{T}$; feasibility indicators $A^c_{ik},A^s_{jk},A_{ijk,t}$; capacities $\mathrm{Cap}_k^t$; parameters $q_{ijk},\bar{E}_{ijk},d_{ik},d_{jk},a_{ik},a_{jk},\alpha,\beta$; repair parameters $\lambda,\mu$}
\KwOut{Pareto-optimal matching solutions $\mathcal{F}^{\star}$}

Use the precomputed feasible session universe $\Omega$\;


$\mathcal{P}_{0}^{init} \leftarrow$ initial population of $N$ individuals using \textbf{SelectSession} from Algorithm~\ref{alg:select-session}\;

$\mathcal{P}_{0} \leftarrow$ repair each individual in $\mathcal{P}_{0}^{init}$ using \textbf{RepairFeasibility} from Algorithm~\ref{alg:repair-feasibility}\;

\For{$g \leftarrow 1$ \KwTo $G$}{
    Evaluate fitness of all individuals in $\mathcal{P}_{g-1}$ using \textbf{EvaluateFitness}$(\mathcal{P}_{g-1},\bar{E},d,a,q,\mathrm{Cap},A^c,A^s,A_{ijk,t},\alpha,\beta)$ from Algorithm~\ref{alg:evaluate-fitness}\;
    
    
    Select parents from $\mathcal{P}_{g-1}$ using tournament selection\;
    
    $\mathcal{P}_{off} \leftarrow$ offspring population after applying session-level crossover followed by mutation with probability $p_m$\;
    
    $\mathcal{P}_{off}^{rep} \leftarrow \emptyset$\;
    
    \ForEach{offspring solution $x \in \mathcal{P}_{off}$}{
        $x^{rep} \leftarrow \textbf{RepairFeasibility}(x,\Omega,\mathrm{Cap},q,\bar{E},d,a,\lambda,\mu)$ from Algorithm~\ref{alg:repair-feasibility}\;
        
        $\mathcal{P}_{off}^{rep} \leftarrow \mathcal{P}_{off}^{rep} \cup \{x^{rep}\}$\;
    }
    
    $\mathcal{P}_{off} \leftarrow \mathcal{P}_{off}^{rep}$\;
    
    Evaluate fitness of all individuals in $\mathcal{P}_{off}$ using \textbf{EvaluateFitness}$(\mathcal{P}_{off},\bar{E},d,a,q,\mathrm{Cap},A^c,A^s,A_{ijk,t},\alpha,\beta)$ from Algorithm~\ref{alg:evaluate-fitness}\;
    
    $\mathcal{R}_{g} \leftarrow \mathcal{P}_{g-1} \cup \mathcal{P}_{off}$\;

$\mathcal{P}_{g} \leftarrow$ select the best $N$ individuals from $\mathcal{R}_{g}$ using non-dominated sorting and crowding distance\;

    

}


$\mathcal{F}^{\star} \leftarrow$ first non-dominated front of $\mathcal{P}_{G}$\;

\Return $\mathcal{F}^{\star}$\;
\end{algorithm*}

\paragraph{Generational Evolution}
The evolutionary search proceeds for $G$ generations. In each generation, the current population $\mathcal{P}_{g-1}$ is evaluated using \textbf{EvaluateFitness} in line~5. The evaluation computes the four objective values: matching coverage $f_1$, transferred energy $f_2$, mobility cost $f_3$, and suitability $f_4$, together with any constraint-violation score. 

Parent solutions are selected in line~6 using binary tournament selection. In each tournament, two solutions are randomly sampled. The solution with the lower non-dominated rank is selected; if both have the same rank, the one with the higher crowding distance is selected. This favours high-quality solutions while preserving diversity.

\paragraph{Variation and Repair}
The offspring population $\mathcal{P}_{off}$ is generated in line~7 using session-level crossover followed by mutation with probability $p_m$. Since crossover and mutation may create infeasible combinations, each offspring solution is repaired in lines~9--12. The repaired solutions are stored in $\mathcal{P}_{off}^{rep}$, and the offspring population is replaced by this repaired population in line~12. The repaired offspring are then evaluated using \textbf{EvaluateFitness} in line~13.


\paragraph{Environmental Selection}
The parent and offspring populations are merged into $\mathcal{R}_g$ in line~14. The next generation $\mathcal{P}_g$ is then formed in line~15 by selecting the best $N$ individuals from $\mathcal{R}_g$ using non-dominated sorting and crowding distance. This selection prioritizes individuals with lower non-dominated rank, while crowding distance is used to preserve diversity among solutions within the same front.

\paragraph{Dominance and Tie-Breaking}
An individual dominates another if it is no worse in all four objectives and strictly better in at least one. Specifically, dominance prefers higher matching coverage $f_1$, higher transferred energy $f_2$, lower mobility cost $f_3$, and higher suitability $f_4$. If two solutions have the same rank and crowding distance, the preference order is higher $f_1$, then higher $f_2$, then lower $f_3$, and finally higher $f_4$.

\paragraph{Final Pareto-Front Output}
\label{pareto-analysis}
After completing $G$ generations, the first non-dominated front of the final population $\mathcal{P}_G$ is extracted in line~16 and assigned to $\mathcal{F}^{\star}$. The algorithm then returns $\mathcal{F}^{\star}$ in line~17. This output provides a set of Pareto-efficient EV--EV matching solutions rather than a single maximum solution. Each solution in $\mathcal{F}^{\star}$ represents a different trade-off among matching coverage, transferred energy, mobility cost, and charging-node suitability, allowing operators to select a preferred solution according to operational priorities.

For direct comparison with single-solution baselines, one balanced NSGA-II solution is selected from $\mathcal{F}^{\star}$ as a post-processing step. Since the objectives have different scales and directions, the selection uses normalized objective values and treats mobility cost as a lower-is-better objective. The balanced solution is chosen to provide the strongest overall compromise across the four objectives, with higher priority given to matching coverage and transferred energy, followed by mobility cost and suitability. This selected solution is reported in the aggregate comparison tables. Separately, objective-specific Pareto representatives are selected from $\mathcal{F}^{\star}$ by choosing the solutions with maximum $f_1$, maximum $f_2$, minimum $f_3$, and maximum $f_4$, respectively.

\subsubsection{\textbf{Feasible Session Selection}}
\label{subsec:pg-nsga2-select-desc}

Algorithm~\ref{alg:select-session} presents the \textbf{SelectSession} procedure, which selects at most one feasible charging session for a given consumer EV decision event $i$ under the current partial solution $x$. The procedure is used during population construction, and it can also support mutation when a new feasible session needs to be selected. A returned session is represented as $s_i^{\star}=(j,k,t)$, which specifies the selected provider $j$, WCP $k$, and time slot $t$ for consumer $i$. If no feasible session exists, the procedure returns $\emptyset$.

\paragraph{Construction strategy}
The procedure constructs sessions in a consumer-centric manner. Although the MILP formulation treats consumers and providers symmetrically, this sequential construction is computationally convenient because each consumer EV decision event is considered once under the current partial solution. Provider-side feasibility is still enforced through the provider reachability indicator $A^s_{jk}$, the temporal feasibility indicator $A_{ijk,t}$, and the condition that provider $j$ must not already be matched in $x$. Any remaining conflicts are handled later by Algorithm~\ref{alg:repair-feasibility}. Therefore, the procedure supports efficient candidate construction while preserving the feasibility structure of the original matching formulation.

\paragraph{Input and Output}
The procedure receives a consumer $i$, the current partial solution $x$, and the sets of providers $\mathcal{S}$, WCPs $\mathcal{K}$, and time slots $\mathcal{T}$. It also receives the feasibility indicators $A^c_{ik}$, $A^s_{jk}$, and $A_{ijk,t}$, WCP-slot capacities $\mathrm{Cap}_k^t$, suitability scores $q_{ijk}$, transferable energy caps $\bar{E}_{ijk}$, the normalization constant $\bar{E}_{\max}$, the scoring mode $m$, and the tie-noise scale $\sigma_{\epsilon}$. The output is the selected feasible session $s_i^{\star}$.

\begin{algorithm}[t]
\small
\caption{Feasible Session Selection}
\label{alg:select-session}
\DontPrintSemicolon


\KwIn{consumer $i$; partial solution $x$; sets $\mathcal{S},\mathcal{K},\mathcal{T}$; feasibility indicators $A^c_{ik},A^s_{jk},A_{ijk,t}$; capacities $\mathrm{Cap}_k^t$; suitability scores $q_{ijk}$; energy caps $\bar{E}_{ijk}$; normalization $\bar{E}_{\max}$; scoring mode $m\in\{\text{suitability},\text{energy},\text{hybrid}\}$; tie-noise scale $\sigma_{\epsilon}$}

\KwOut{Selected feasible session $s_i^{\star}$}

Initialize $s_i^{\star} \leftarrow\emptyset$\;

Initialize feasible candidate set $\mathcal{F}_{i}\leftarrow\emptyset$\;

\If{\textnormal{consumer} $i$ \textnormal{is already matched in} $x$}{
    \Return $\emptyset$\;
}

                

                

\ForEach{\textnormal{WCP } $k\in\mathcal{K}$ \textnormal{ and } $A^c_{ik}=1$}{
    \ForEach{\textnormal{provider } $j\in\mathcal{S}$ \textnormal{ and } $A^s_{jk}=1$}{
        \If{\textnormal{provider } $j$ \textnormal{ is not already matched in } $x$}{
            \ForEach{\textnormal{time slot } $t\in\mathcal{T}$ \textnormal{ and } $A_{ijk,t}=1$}{
                $u_{kt}\leftarrow$ \textnormal{number of sessions in } $x$ \textnormal{ assigned to WCP--slot pair } $(k,t)$\;
                
                \If{$u_{kt}<\mathrm{Cap}_{k}^{t}$}{
                    \textnormal{Add session } $(j,k,t)$ \textnormal{ to } $\mathcal{F}_{i}$\;
                }
            }
        }
    }
}

\If{$\mathcal{F}_{i}=\emptyset$}{
    \Return $\emptyset$\;
}

\ForEach{session $s \in\mathcal{F}_{i}$}{
    \uIf{$m=\mathrm{suitability}$}{
        $\epsilon \sim U(0,\sigma_{\epsilon})$\;
        $\pi(s)\leftarrow q_{ijk}+\epsilon$\;
    }
    \uElseIf{$m=\mathrm{energy}$}{
        $\epsilon \sim U(0,\sigma_{\epsilon})$\;
        $\pi(s)\leftarrow \bar{E}_{ijk}+\epsilon$\;
    }
    \Else{
        $\epsilon \sim U(0,\sigma_{\epsilon})$\;
        $\pi(s)\leftarrow 0.4q_{ijk}+0.6\bar{E}_{ijk}/\bar{E}_{\max}+\epsilon$\;
    }
}

$s_i^{\star}\leftarrow \arg\max_{\forall s\in\mathcal{F}_{i}}\pi(s)$\;


\Return $s_i^{\star}$\;
\end{algorithm}

\paragraph{Feasible candidate set initialization (Line 2)}
The feasible candidate set $\mathcal{F}_i$ is initialized as empty in line~2. This set stores all feasible provider--WCP--time combinations $(j,k,t)$ that can be assigned to consumer $i$ under the current partial solution.

\paragraph{Consumer availability check (Lines 3--4)}
The procedure first checks whether consumer $i$ is already matched in the current partial solution $x$ (line~3). If consumer $i$ already appears in an existing session, the procedure returns $\emptyset$ (line~4). This prevents the same consumer from being assigned to more than one session during incremental solution construction.

\paragraph{Feasible session enumeration (Lines 5--11)}
The procedure enumerates candidate sessions by checking WCPs, providers, and time slots. For each WCP $k$ that is reachable by consumer $i$ ($A^c_{ik}=1$), the algorithm considers each provider $j$ that can reach the same WCP ($A^s_{jk}=1$) in lines~5--6. It then checks whether provider $j$ is not already matched in the partial solution $x$ (line~7). For each temporally feasible time slot $t$ satisfying $A_{ijk,t}=1$ (line~8), the algorithm computes the current WCP-slot usage $u_{kt}$ from the sessions already assigned in $x$ (line~9). If the current usage is below the available capacity $\mathrm{Cap}_k^t$, the session $(j,k,t)$ is added to $\mathcal{F}_i$ (lines~10--11). Hence, every session added to $\mathcal{F}_i$ satisfies consumer reachability, provider reachability, temporal feasibility, provider exclusivity, and WCP-slot capacity availability.

\paragraph{Empty candidate handling (Lines 12--13)}
After candidate enumeration, the procedure checks whether $\mathcal{F}_i$ is empty (line~12). If no feasible session has been found for consumer $i$, the procedure returns $\emptyset$ (line~13). In the main NSGA-II algorithm, this means consumer $i$ remains unmatched in the current individual unless later variation operators create a feasible assignment.

\paragraph{Candidate scoring (Lines 14--23)}
Each feasible session $s=(j,k,t)\in\mathcal{F}_i$ is assigned a utility value $\pi(s)$ in lines~14--23. If the scoring mode is \textit{suitability}, the utility is computed as $q_{ijk}+\epsilon$ (lines~15--17), giving priority to sessions with higher prediction-derived suitability. If the scoring mode is \textit{energy}, the utility is computed as $\bar{E}_{ijk}+\epsilon$ (lines~18--20), giving priority to sessions with higher transferable energy. Otherwise, the hybrid mode combines both signals as $0.4q_{ijk}+0.6\bar{E}_{ijk}/\bar{E}_{\max}+\epsilon$ (lines~21--23). The noise term $\epsilon$ is sampled from a uniform distribution and is used to break ties and introduce diversity among candidate sessions.

\paragraph{Session selection and return value (Lines 24--25)}
After all feasible sessions have been scored, the procedure selects the session with the highest utility value in line~24:
\[
s_i^{\star} \leftarrow \arg\max_{s\in\mathcal{F}_i}\pi(s).
\]
The selected session $s_i^{\star}$ is returned in line~25. Since $s_i^{\star}$ is selected from $\mathcal{F}_i$, it is feasible with respect to the session-level checks performed during enumeration. When used in the main NSGA-II algorithm, this returned session is combined with consumer $i$ to form the full assignment $(i,j,k,t)$.

\subsubsection{\textbf{Feasibility Repair}} \label{subsec:pg-nsga2-repair-desc} Algorithm~\ref{alg:repair-feasibility} presents the \textbf{RepairFeasibility} procedure, which converts a candidate session set $x$ into a repaired EV--EV matching solution. The procedure is applied after population initialization and after session-level crossover or mutation, because these operations may generate duplicated sessions, infeasible tuples, multiple assignments for the same consumer or provider, or WCP--slot capacity violations. Before applying feasibility repair, each session is assigned a deterministic repair score: \begin{equation} S(i,j,k,t)= q_{ijk} \cdot \bar{E}_{ijk} -\lambda\frac{d_{ik}+d_{jk}}{D_{\max}} -\mu\frac{|a_{ik}-a_{jk}|}{W_{\max}} . \label{eq:repair-score} \end{equation} The first term favours sessions with high prediction-derived suitability and high transferable energy, while the second and third terms penalize longer travel distance and larger arrival-time mismatch, respectively. The weights $\lambda$ and $\mu$ control the relative importance of distance and temporal-mismatch penalties during repair.
For compact notation in Algorithm~\ref{alg:repair-feasibility}, we use $S(s)$ to denote $S(i,j,k,t)$ for a session $s=(i,j,k,t)$.

\begin{algorithm}[t]
\small
\caption{Feasibility Repair for EV--EV Matching}
\label{alg:repair-feasibility}
\DontPrintSemicolon

\KwIn{candidate solution $x$ (set of sessions $(i,j,k,t)$); sets $\mathcal{C},\mathcal{S},\mathcal{K},\mathcal{T}$; feasible session universe $\Omega$; capacities $\mathrm{Cap}_k^t$}

\KwOut{Repaired feasible solution $x^{\mathrm{rep}}$}

Initialize $x^{\mathrm{rep}}\leftarrow\emptyset$\;

Define repair score $S(i,j,k,t)$ using Eq.~\eqref{eq:repair-score}\;


Update $x$ by removing duplicated sessions and sessions not contained in $\Omega$\;


\ForEach{consumer $i \in \mathcal{C}$}{
    $\mathcal{X}_{i} \leftarrow$ Identify all sessions $s$ in $x$ involving consumer $i$\;
    
    \If{$|\mathcal{X}_{i}|>1$}{
        $s_i^{\star} \leftarrow \arg\max_{s\in\mathcal{X}_{i}} S(s)$\;




        $x \leftarrow (x\setminus\mathcal{X}_{i}) \cup \{s_i^{\star}\}$\;

    }
}

\ForEach{provider $j \in \mathcal{S}$}{
$\mathcal{X}_{j} \leftarrow$ Identify all sessions $s$ in $x$ involving provider $j$\;

    \If{$|\mathcal{X}_{j}|>1$}{
        $s_j^{\star} \leftarrow \arg\max_{s\in\mathcal{X}_{j}} S(s)$\;
        $x \leftarrow (x\setminus\mathcal{X}_{j}) \cup \{s_j^{\star}\}$\;

    }
}

    

\ForEach{WCP--slot pair $(k,t)\in\mathcal{K}\times\mathcal{T}$}{

    $\mathcal{X}_{kt} \leftarrow$ Identify all sessions $s$ in $x$ assigned to $(k,t)$\;
    
    \If{$|\mathcal{X}_{kt}|>\mathrm{Cap}_{k}^{t}$}{
        Sort $\mathcal{X}_{kt}$ in descending order of $S(s)$\;

        
        $\mathcal{X}_{kt}^{\mathrm{keep}} \leftarrow$ Keep the highest-scoring $\mathrm{Cap}_{k}^{t}$ sessions\;



    $x \leftarrow (x\setminus\mathcal{X}_{kt}) \cup \mathcal{X}_{kt}^{\mathrm{keep}}$\;
    }
}

$x^{\mathrm{rep}} \leftarrow x$\;
\Return $x^{\mathrm{rep}}$\;
\end{algorithm}

\paragraph{Input and Output}
The procedure receives a candidate solution $x$, represented as a set of sessions $(i,j,k,t)$; the consumer, provider, WCP, and time-slot sets $\mathcal{C}$, $\mathcal{S}$, $\mathcal{K}$, and $\mathcal{T}$; the feasible session universe $\Omega$; and the WCP--slot capacities $\mathrm{Cap}_k^t$. The score $S(i,j,k,t)$ used for conflict resolution is defined in Eq.~\eqref{eq:repair-score}. The output is a repaired feasible solution $x^{\mathrm{rep}}$.

\paragraph{Repair score definition (Line 2)}
The algorithm first uses the repair score $S(i,j,k,t)$ defined in Eq.~\eqref{eq:repair-score}. This score ranks competing sessions during repair by favouring sessions with higher suitability and transferable energy, while penalising larger travel distance and arrival-time mismatch.


\paragraph{Duplicate and infeasible-session removal (Line 3)}
After defining the repair score, the procedure updates the candidate solution $x$ in line~3 by removing duplicated sessions and sessions that are not contained in the feasible session universe $\Omega$. Since $\Omega$ contains only spatially and temporally feasible sessions, this step removes assignments that violate the basic feasibility requirements before one-to-one and capacity conflicts are handled.



\paragraph{Consumer one-to-one repair (Lines 4--8)}
For each consumer $i\in\mathcal{C}$, the algorithm identifies the subset $\mathcal{X}_i$ of sessions in $x$ involving consumer $i$ in line~5. If more than one session involves the same consumer, the session with the highest repair score is selected as $s_i^{\star}$ in line~7. The algorithm then updates $x$ in line~8 by replacing all sessions in $\mathcal{X}_i$ with the selected session $s_i^{\star}$, i.e., $x \leftarrow (x\setminus\mathcal{X}_i)\cup\{s_i^{\star}\}$. After this step, each consumer EV decision event appears in at most one session.



\paragraph{Provider one-to-one repair (Lines 9--13)}
The provider-side repair follows the same logic as the consumer-side repair. For each provider $j\in\mathcal{S}$, the algorithm identifies the subset $\mathcal{X}_j$ of sessions in $x$ involving provider $j$ in line~10. If multiple sessions involve the same provider, the highest-scoring session is selected as $s_j^{\star}$ in line~12. The algorithm then updates $x$ in line~13 by replacing all sessions in $\mathcal{X}_j$ with the selected session $s_j^{\star}$, i.e., $x \leftarrow (x\setminus\mathcal{X}_j)\cup\{s_j^{\star}\}$. This ensures that each provider EV decision event appears in at most one selected session.

\paragraph{WCP--slot capacity repair (Lines 14--19)}
The algorithm then enforces the WCP--slot capacity constraint. For each WCP--slot pair $(k,t)\in\mathcal{K}\times\mathcal{T}$, it identifies the set $\mathcal{X}_{kt}$ of sessions assigned to that pair in line~15. If the number of assigned sessions exceeds $\mathrm{Cap}_k^t$, the sessions are sorted in descending order of the repair score in line~17. The highest-scoring $\mathrm{Cap}_k^t$ sessions are retained in $\mathcal{X}_{kt}^{\mathrm{keep}}$ in line~18. The algorithm then updates $x$ in line~19 by replacing all sessions in $\mathcal{X}_{kt}$ with the retained set $\mathcal{X}_{kt}^{\mathrm{keep}}$, i.e., $x \leftarrow (x\setminus\mathcal{X}_{kt})\cup\mathcal{X}_{kt}^{\mathrm{keep}}$. This keeps the WCP--slot usage within its capacity while retaining the strongest sessions according to the repair score.

\paragraph{Return value (Lines 20--21)}
After duplicate removal, infeasible-session filtering, one-to-one repair, and WCP--slot capacity repair, the remaining session set is assigned to $x^{\mathrm{rep}}$ in line~20. The algorithm then returns $x^{\mathrm{rep}}$ as the repaired solution in line~21.

\subsubsection{\textbf{Constraint-Aware Fitness Evaluation}}
\label{subsec:pg-nsga2-fit-desc}

Algorithm~\ref{alg:evaluate-fitness} presents the \textbf{EvaluateFitness} procedure, which evaluates each candidate EV--EV matching solution in a population. For every solution $x\in\mathcal{P}$, the algorithm computes the four objective values: matching coverage $f_1(x)$, transferred energy $f_2(x)$, mobility cost $f_3(x)$, and prediction-derived suitability $f_4(x)$. It also computes an aggregate constraint-violation score $\mathrm{Viol}(x)$. The resulting fitness record $\mathcal{F}$ is used by NSGA-II for non-dominated sorting, crowding-distance assignment, and environmental selection.

\paragraph{Input and Output}
The procedure takes as input a population $\mathcal{P}$ of candidate matching solutions; the sets of consumers $\mathcal{C}$, providers $\mathcal{S}$, WCPs $\mathcal{K}$, and time slots $\mathcal{T}$; transferable-energy caps $\bar{E}_{ijk}$; suitability scores $q_{ijk}$; travel distances $d_{ik}$ and $d_{jk}$; arrival times $a_{ik}$ and $a_{jk}$; normalization constants $D_{\max}$ and $W_{\max}$; objective weights $\alpha$ and $\beta$; WCP--slot capacities $\mathrm{Cap}_k^t$; and feasibility indicators $A^c_{ik}$, $A^s_{jk}$, and $A_{ijk,t}$. The output is a fitness record $\mathcal{F}$ containing, for each solution $x$, the tuple $(x,f_1(x),f_2(x),f_3(x),f_4(x),\mathrm{Viol}(x))$.

\begin{algorithm}[t]
\small
\caption{Constraint-Aware Fitness Evaluation}
\label{alg:evaluate-fitness}
\DontPrintSemicolon
\KwIn{population $\mathcal{P}$; sets $\mathcal{C},\mathcal{S},\mathcal{K},\mathcal{T}$; energy caps $\bar{E}_{ijk}$; suitability $q_{ijk}$; distances $d_{ik},d_{jk}$; arrival times $a_{ik},a_{jk}$; normalization constants $D_{\max},W_{\max}$; weights $\alpha,\beta$; capacities $\mathrm{Cap}_k^t$; feasibility indicators $A^c_{ik},A^s_{jk},A_{ijk,t}$}
\KwOut{Fitness record $F$}


Initialize fitness record $F\leftarrow\emptyset$\;

\ForEach{solution $x \in \mathcal{P}$}{
    Initialize $f_1(x)\leftarrow 0$, $f_2(x)\leftarrow 0$, $f_3(x)\leftarrow 0$, and $f_4(x)\leftarrow 0$\;
    
    Initialize violation counters $v_{C1},v_{C2},v_{C3},v_{C4},v_{C5},v_{C6}\leftarrow 0$\;
    
    \ForEach{session $(i,j,k,t)\in x$}{
        $f_1(x)\leftarrow f_1(x)+1$\;
        $f_2(x)\leftarrow f_2(x)+\bar{E}_{ijk}$\;
        $f_3(x)\leftarrow f_3(x)+
        \left[
        \alpha \frac{d_{ik}+d_{jk}}{D_{\max}}
        +
        \beta \frac{|a_{ik}-a_{jk}|}{W_{\max}}
        \right]$\;
        $f_4(x)\leftarrow f_4(x)+q_{ijk}$\;
    }
    
    Construct $x^{\mathrm{uniq}}$ by retaining one copy of each distinct session in $x$\;
    $v_{C6}\leftarrow |x|-|x^{\mathrm{uniq}}|$\;
    
    \ForEach{session $(i,j,k,t)\in x^{\mathrm{uniq}}$}{
        \If{$A^c_{ik}=0$ or $A^s_{jk}=0$}{
            $v_{C4}\leftarrow v_{C4}+1$\;
        }
        \If{$A_{ijk,t}=0$}{
            $v_{C5}\leftarrow v_{C5}+1$\;
        }
    }
    
    \ForEach{consumer $i\in\mathcal{C}$}{
        $n_i\leftarrow$ number of sessions in $x^{\mathrm{uniq}}$ involving consumer $i$\;
        \If{$n_i>1$}{
            $v_{C1}\leftarrow v_{C1}+(n_i-1)$\;
        }
    }
    
    \ForEach{provider $j\in\mathcal{S}$}{
        $n_j\leftarrow$ number of sessions in $x^{\mathrm{uniq}}$ involving provider $j$\;
        \If{$n_j>1$}{
            $v_{C2}\leftarrow v_{C2}+(n_j-1)$\;
        }
    }
    
    \ForEach{WCP--slot pair $(k,t)\in\mathcal{K}\times\mathcal{T}$}{
        $n_{kt}\leftarrow$ number of sessions in $x^{\mathrm{uniq}}$ assigned to $(k,t)$\;
        \If{$n_{kt}>\mathrm{Cap}_k^t$}{
            $v_{C3}\leftarrow v_{C3}+(n_{kt}-\mathrm{Cap}_k^t)$\;
        }
    }
    
    $\mathrm{Viol}(x)\leftarrow v_{C1}+v_{C2}+v_{C3}+v_{C4}+v_{C5}+v_{C6}$\;
    
    Add $(x,f_1(x),f_2(x),f_3(x),f_4(x),\mathrm{Viol}(x))$ to $F$\;
}

\Return $F$\;
\end{algorithm}

\paragraph{Fitness record initialization (Line 1)}
The algorithm first initializes the fitness record $F$ as an empty set in line~1. This record stores the objective values and violation score of every candidate solution evaluated in the population.

\paragraph{Per-solution initialization (Lines 2--4)}
For each candidate solution $x\in\mathcal{P}$, the algorithm initializes the four objective values $f_1(x)$, $f_2(x)$, $f_3(x)$, and $f_4(x)$ to zero in line~3. It also initializes the constraint-violation counters $v_{C1}$, $v_{C2}$, $v_{C3}$, $v_{C4}$, $v_{C5}$, and $v_{C6}$ to zero in line~4. These counters correspond to consumer decision-event uniqueness, provider decision-event uniqueness, WCP--slot capacity, spatial feasibility, temporal feasibility, and duplicate-session violations, respectively.

\paragraph{Objective computation (Lines 5--9)}
The algorithm then iterates over each selected session $(i,j,k,t)\in x$ in line~5. For every session, the matching objective $f_1(x)$ is increased by one in line~6. The transferred-energy objective $f_2(x)$ accumulates the transferable energy $\bar{E}_{ijk}$ in line~7. The mobility-cost objective $f_3(x)$ accumulates the weighted and normalized distance and arrival-time mismatch in line~8. The suitability objective $f_4(x)$ accumulates the prediction-derived suitability score $q_{ijk}$ in line~9. Together, these computations evaluate the four objectives used by the NSGA-II search.

\paragraph{Duplicate-session violation (Lines 10--11)}
The algorithm creates a duplicate-free version of the solution, denoted by $x^{\mathrm{uniq}}$, in line~10. The duplicate-session violation counter $v_{C6}$ is then computed in line~11 as the difference between the size of the original solution $x$ and the size of $x^{\mathrm{uniq}}$. A positive value indicates that duplicated session tuples exist in the candidate solution.

\paragraph{Spatial and temporal feasibility violations (Lines 12--16)}
For each unique session $(i,j,k,t)\in x^{\mathrm{uniq}}$, the algorithm checks spatial and temporal feasibility. If either the consumer-side spatial feasibility indicator $A^c_{ik}$ or the provider-side spatial feasibility indicator $A^s_{jk}$ is zero, the spatial-feasibility violation counter $v_{C4}$ is increased in lines~13--14. If the temporal-feasibility indicator $A_{ijk,t}$ is zero, the temporal-feasibility violation counter $v_{C5}$ is increased in lines~15--16.

\paragraph{Consumer one-to-one violation (Lines 17--20)}
The algorithm checks the consumer one-to-one constraint in lines~17--20. For each consumer $i\in\mathcal{C}$, it counts the number of sessions in $x^{\mathrm{uniq}}$ involving that consumer. If consumer $i$ appears in more than one session, the excess number of assignments, $n_i-1$, is added to $v_{C1}$. This quantifies the extent to which the candidate solution violates consumer decision-event uniqueness.

\paragraph{Provider one-to-one violation (Lines 21--24)}
The provider one-to-one constraint is checked in the same way in lines~21--24. For each provider $j\in\mathcal{S}$, the algorithm counts the number of sessions involving that provider. If provider $j$ appears in more than one session, the excess number of assignments, $n_j-1$, is added to $v_{C2}$.

\paragraph{WCP--slot capacity violation (Lines 25--28)}
The algorithm checks WCP--slot capacity constraints in lines~25--28. For each WCP--slot pair $(k,t)\in\mathcal{K}\times\mathcal{T}$, it counts the number of sessions assigned to that pair. If the count $n_{kt}$ exceeds the available capacity $\mathrm{Cap}_k^t$, the excess number of sessions, $n_{kt}-\mathrm{Cap}_k^t$, is added to $v_{C3}$.

\paragraph{Violation aggregation and fitness storage (Lines 29--30)}
The aggregate violation score is computed in line~29 as
\[
\mathrm{Viol}(x)=v_{C1}+v_{C2}+v_{C3}+v_{C4}+v_{C5}+v_{C6}.
\]
The algorithm then stores the complete fitness tuple $(x,f_1(x),f_2(x),f_3(x),f_4(x),\mathrm{Viol}(x))$ in the fitness record $F$ in line~30.

\paragraph{Return value (Line 31)}
After all solutions in $\mathcal{P}$ have been evaluated, the algorithm returns the fitness record $F$ in line~31. This record provides the objective and feasibility information required for NSGA-II ranking, diversity preservation, and selection.

In the proposed implementation, candidate solutions are repaired before fitness evaluation so that infeasible assignments are adjusted with respect to the feasible opportunity set, one-to-one matching, temporal feasibility, and WCP capacity. Therefore, evaluated individuals are expected to be feasible after repair. The violation score computed during constraint-aware fitness evaluation is retained as a diagnostic safeguard to detect any residual infeasibility that may remain after repair, rather than as a primary constrained-dominance ranking criterion.

\subsection{Baseline Methods}
\label{sec:baselines}

All baselines operate on the same feasible matching set induced by the constraints in Section~\ref{sec:cons}; they differ only in their match selection strategies.



\subsubsection{Random Matching}
Random Matching samples feasible tuples uniformly at random and adds them sequentially while maintaining consumer decision-event uniqueness, provider decision-event uniqueness, and WCP--slot capacity feasibility until no further feasible tuple can be added. To reduce the variance caused by a single random rollout, Random Matching is executed for 30 independent rollouts within each evaluation seed. The best feasible solution from these 30 rollouts is retained for that seed using the same reporting criterion adopted for single-solution comparison. The final Random results are then reported as \emph{mean $\pm$ standard deviation} across the five retained seed-level solutions.

\subsubsection{Greedy Matching}
Greedy Matching iteratively selects the feasible tuple maximizing a scalar utility derived directly from the MILP objectives:
\begin{equation}
\begin{aligned}
\pi(i,j,k,t) =\;& w_E\, \bar{E}_{ijk} + w_q\, q_{ijk} \\
&- w_m \left[ \alpha \frac{d_{ik}+d_{jk}}{D_{\max}} + \beta \frac{|a_{ik}-a_{jk}|}{W_{\max}} \right]
\end{aligned}
\label{eq:greedy_priority}
\end{equation}
where $\bar{E}_{ijk}$ is the transferable energy (Eq.~\ref{eq:transfer_cap}), the bracketed term is the normalized mobility cost from $f_3$ (Eq.~\ref{eq:obj_mobility}) with coefficients $\alpha = 0.7$ and $\beta = 0.3$, $q_{ijk}$ is the predicted pairwise suitability, and $w_E, w_m, w_q \ge 0$ are weights that control the relative emphasis on the three terms. After each selection, all conflicting tuples (violating consumer or provider exclusivity or slot capacity) are removed, and the process repeats until no feasible tuples remain.


\subsubsection{Proximal Policy Optimization (PPO)}
PPO trains an actor--critic deep reinforcement learning agent to sequentially construct a feasible EV--EV matching. At each decision step, the agent observes the current environment state, including unmatched consumers, unmatched providers, and remaining feasible sessions, and selects a feasible consumer--provider--WCP--time assignment according to a stochastic policy network. After each selected assignment, the environment updates the remaining feasible matching space by enforcing the one-to-one matching, spatial feasibility, temporal feasibility, WCP-slot capacity, and charging-rate constraints. The policy and value networks are optimized using the clipped surrogate objective proposed by Schulman et al.~\cite{schulman2017proximal}, with generalized advantage estimation (GAE) used to reduce variance in the policy-gradient update~\cite{schulman2015high}. The agent is trained to maximize a scalarized reward that combines the four objectives ($f_1$--$f_4$) using fixed weights, subject to all feasibility and capacity constraints. The final matching is constructed by rolling out the trained policy until no feasible actions remain.


\subsection{State-of-the-art Methods}
\label{sec:soa}
\subsubsection{Double Auction Matching Algorithm (DAMA)}

The DAMA \cite{C17} performs bid--ask clearing independently within each wireless charging pad and time slot $(k,t)$, subject to the same feasibility and capacity constraints defined in Section~\ref{sec:cons}. 
Consumers submit bids that reflect their charging demand and mobility cost, while providers submit asks that reflect their energy supply and mobility cost. 
Within each $(k,t)$, consumer bids are sorted in descending order and provider asks in ascending order, and matches are formed sequentially as long as bid values meet or exceed ask values, while enforcing one-to-one matching and slot capacity constraints. 
The DAMA mechanism thus focuses on clearing feasible local trades and maximizing auction surplus at each $(k,t)$, rather than optimizing global Pareto-optimal trade-offs across objectives.

\subsubsection{Bichromatic Mutual Nearest Neighbor (BMNN)}

The BMNN baseline implements the original distance-based matching algorithm proposed by Yucel et al.~\cite{P-T-129}. In this approach, each unmatched consumer identifies its nearest feasible provider (minimizing the sum of travel distances), and each unmatched provider identifies its nearest feasible consumer. A match is formed only if both parties mutually select each other as their nearest neighbor, and the match does not violate any one-to-one matching or slot capacity constraints. This process is repeated iteratively until no further mutual nearest neighbor pairs can be found. BMNN thus prioritizes spatial proximity in matching, without considering price or energy maximization objectives.




\section{Experiment and Result Analysis}
\label{sec:results}

\subsection{Experimental Setup}
\label{sec:setup}

We evaluate the proposed framework using predicted EV--EV trading data for the Chicago metropolitan area on January~1, 2016 (Friday), generated by EVNextTrade~\cite{rahman2025evnexttrade}. The dataset contains 431 consumer EV decision events and 388 provider EV decision events generated by 204 and 210 distinct physical EVs, respectively, distributed across 128 charging stations over a 24-hour period. Each EV decision event is associated with up to 31 ranked candidate charging stations based on proximity. For computational tractability, we sample 150 consumer EV decision events and 150 provider EV decision events for optimization, corresponding to 300 sampled EV decision events generated by 186 unique physical EVs. Table~\ref{tab:experimental_setup} summarizes the key parameter settings.

Matching coverage is computed at the EV decision-event level. Unique physical-EV counts are reported only to describe the diversity of the sampled instance and are not used as the denominator for matching coverage.

For each experimental instance, the feasible opportunity set $\Omega$ is reconstructed to enforce identical spatial and temporal feasibility constraints, including distance limits, waiting tolerances, and charging-pad capacity. This ensures that all algorithms operate on the same feasible matching space.

The hyperparameters of NSGA-II were determined through a one-parameter-at-a-time sensitivity analysis conducted on the same 150-consumer/150-provider EV decision-event instance. In this analysis, one parameter was varied while the remaining parameters were fixed at their default settings. Each configuration was evaluated using three independent random seeds, and the results were assessed using a balanced multi-objective criterion rather than matching percentage alone. The population size and the number of generations were selected from $\{50,100,150,200,250\}$, while the crossover and mutation rates were varied over $\{0.2,0.4,0.6,0.8\}$. Based on the sensitivity analysis results presented in Table~\ref{tab:sensitivity}, the final configuration was set to a crossover rate of $0.2$, a mutation rate of $0.6$, $250$ generations, and a population size of $250$. These settings were selected because they maintain competitive matching performance while supporting transferred energy, mobility cost, suitability, and Pareto-front exploration.

To ensure a fair and rigorous comparison, all algorithms are evaluated using identical objective definitions ($f_1$--$f_4$), the same five feasibility constraints defined in Eqs. (1)--(12), and a common dataset derived from the same input instance. This design ensures that observed performance differences arise from algorithmic search behaviour rather than relaxed constraints or inconsistent problem formulations.

For the main comparison experiments, stochastic algorithms were evaluated over five independent random seeds, $\{1,2,3,4,42\}$, and results are reported as mean $\pm$ standard deviation. Deterministic methods, including BMNN and DAMA, produce identical outputs for the fixed input instance and therefore have zero standard deviation across repeated runs. For each NSGA-II run, the algorithm returns a final non-dominated front containing multiple Pareto-efficient trade-off solutions, as described in Section~\ref{pareto-analysis}. Since the comparison methods return a single solution, one representative balanced solution is selected from this front for direct comparison. This selection is performed as a post-processing step by normalizing the objective values and jointly considering matching coverage, transferred energy, mobility cost, and suitability. The selected solution is the one that provides the strongest overall compromise across these objectives, rather than the best value for any single objective alone. Accordingly, the aggregate comparison tables report the mean $\pm$ standard deviation of the balanced NSGA-II solutions selected across the five independent seeds. In contrast, the Pareto-front tables report representative objective-specific solutions from the final non-dominated front of the seed 42 run, including the best-$f_1$, best-$f_2$, best-$f_3$, and best-$f_4$ solutions, to illustrate the range of trade-offs available to decision-makers.

The evaluation proceeds in two stages. In the \emph{baseline study}, NSGA-II is compared against three reference algorithms of increasing sophistication: (i)~\textbf{Random}, which uniformly samples feasible tuples while preserving feasibility; (ii)~\textbf{Greedy}, which sequentially selects feasible sessions using a scalar utility aligned with the optimization objectives; and (iii)~\textbf{PPO}, a proximal-policy-optimization agent trained with an actor--critic architecture. These baselines span the spectrum from stochastic sampling through single-objective heuristic search to deep reinforcement learning, thereby isolating the contribution of multi-objective evolutionary search.

In the \emph{state-of-the-art study}, NSGA-II is compared against two established algorithms for EV--EV energy trading: (i)~\textbf{DAMA} (Double Auction Matching Algorithm)~\cite{C17}, a market-based double-sided auction mechanism that matches buyers and sellers based on bid--ask clearing; and (ii)~\textbf{BMNN} (Bichromatic Mutual Nearest Neighbor)~\cite{P-T-129}, a proximity-based matching algorithm that pairs consumers and providers by minimizing combined spatial distance. These methods are chosen because they represent two major paradigms in the EV--EV trading literature: auction-based coordination and proximity-based matching.


For both DAMA and BMNN, we implement the state-of-the-art algorithms as faithfully as possible while applying the same feasible opportunity set $\Omega$, feasibility constraints, and objective definitions used for NSGA-II. Because the original BMNN and DAMA formulations do not include all journey-aware feasibility constraints considered in this study, they are adapted to operate on the same feasible opportunity set $\Omega$. This adaptation ensures constraint-level fairness, so that performance differences arise from the matching logic of each algorithm rather than from differences in feasibility filtering or constraint relaxation. At the same time, the original decision logic of each method is preserved: BMNN remains a proximity-based mutual nearest-neighbor matching method, while DAMA remains a bid--ask clearing mechanism for auction-based provider--consumer matching.

\textbf{Travel Time Computation.}
In the formulation, arrival times are expressed using route-level segment travel times. Since segment-level routing coverage in the dataset is incomplete and not all road segments form a fully connected graph, we adopt a hybrid travel-time model. When EV~$i$ and charging station~$k$ are connected through the available segment network, the arrival time is computed using the shortest path over connected road segments. Otherwise, the arrival time is approximated using the direct road-network distance $d_{ik}$ divided by a time-dependent aggregate traffic speed. Let $\pi_{ik}$ denote the shortest connected path from EV~$i$ to charging station~$k$, when such a path exists.


\begin{equation}
a_{ik} =
\begin{cases}
t_i^{\mathrm{dep}} + \displaystyle\sum_{(u,v)\in \pi_{ik}} \frac{\ell_{uv}}{v_{uv}(t_i^{\mathrm{dep}})}, & \text{if connected}, \\[6pt]
t_i^{\mathrm{dep}} + \displaystyle\frac{d_{ik}}{v(t_i^{\mathrm{dep}})}, & \text{otherwise}.
\end{cases}
\label{eq:arrival_impl}
\end{equation}
where $t_i^{\mathrm{dep}}$ is the recorded timestamp of EV~$i$ at its current location, $\ell_{uv}$ denotes the length of road segment $(u,v)$, and $v_{uv}(t)$ denotes the travel speed on that segment at time $t$. For disconnected cases, the fallback term uses the following aggregate speed function:
\begin{equation}
v(t) =
\begin{cases}
\mathcal{U}(20, 30) & \text{peak: } t \in [7,9) \cup [16,19) \\[3pt]
\mathcal{U}(30, 40) & \text{midday: } t \in [10,16) \\[3pt]
\mathcal{U}(35, 50) & \text{off-peak: otherwise}
\end{cases}
\;\text{(km/h)}
\label{eq:traffic_impl}
\end{equation}
where $\mathcal{U}(a,b)$ denotes a uniform random variable over $[a,b]$. Equation~\eqref{eq:traffic_impl} defines a stylized time-of-day traffic model for simulation. It is loosely informed by the temporal variability of Chicago-like urban traffic conditions reported in prior sources \cite{trafic_tomtom, chicago_traffic_tracker_segments, chicago_traffic_tracker_segments_origin}, but it does not aim to reproduce segment-level or region-level traffic measurements exactly. The stochastic component captures day-to-day variation in traffic conditions. This hybrid strategy allows the framework to exploit segment-level routing information when available while remaining robust to incomplete connectivity in the extracted road graph. Provider arrival times $a_{jk}$ are computed analogously.

\begin{table}[t]
\centering
\caption{Experimental Parameter Settings}
\label{tab:experimental_setup}
\scriptsize
\renewcommand{\arraystretch}{1.1}
\setlength{\tabcolsep}{3pt}
\begin{tabular}{p{17mm}p{25mm}rc}
\hline
\textbf{Parameter} & \textbf{Description} & \textbf{Value} & \textbf{Source} \\
\hline
\multicolumn{4}{l}{\textit{Dataset \& Temporal Scope}} \\
Region & Study area & Chicago, USA & Data \\
Date & Recommendation day & Jan 1, 2016 (Fri) & Data \\
$T$ & Time horizon & 24 h (00:00--23:59) & Set \\
$\Delta$ & Time-slot duration & 1 h & Set \\
Data source & Recommended EV sessions & EVNextTrade Study & Data \\
\hline
\multicolumn{4}{l}{\textit{Dataset Scale}} \\
$|\mathcal{C}|$ & Consumer EV decision events  & 431 & Data \\
$|\mathcal{S}|$ & Provider EV decision events & 388 & Data \\
$|\mathcal{K}|$ & WCPs & 128 & Data \\
EV models & Distinct vehicle types & 9 & Data \\
Candidates & Ranked stations per EV decision event & [3, 31], $\mu{\approx}10$ & Data \\
\hline
{\textit{Sampled Instance}} \\
$N_C$ & Sampled consumer EV decision events & 150 & Set \\
$N_S$ & Sampled provider EV decision events & 150 & Set \\
$N_C + N_S$ & Total sampled EV decision events & 300 & Set \\
Unique consumer-side EVs & Distinct physical EVs generating sampled consumer events & 109 & Data \\
Unique provider-side EVs & Distinct physical EVs generating sampled provider events & 122 & Data \\
Unique EVs & Distinct physical EVs in sampled instance & 186 & Data \\
\hline
\multicolumn{4}{l}{\textit{Energy Parameters}} \\
$E^{\text{dem}}_c$ & Energy demand & [10.52, 96.6] kWh, $\mu{=}40.50$ & Data \\
$E^{\text{sup}}_s$ & Energy supply & [0.0, 76.26] kWh, $\mu{=}29.33$ & Data \\
$B$ & Battery capacity & [30, 118] kWh, $\mu{=}67.6$ & Data \\
$P_k$ & Charging speed & [3.7, 350] kW, $\mu{=}33.0$ & Data \\
\hline
\multicolumn{4}{l}{\textit{Spatial Parameters}} \\
$d_{i,k}$ & Distance to station & [0.05, 6.85] km, $\mu{=}1.1$ & Data \\
$D_{\max}$ & Max matching distance & 10 km & Set \\
$D^C_{\max}$, $D^S_{\max}$ & Spatial radius (C/S) & 10 km & Set \\
\hline
\multicolumn{4}{l}{\textit{Traffic Parameters}} \\
$v_{\text{peak}}$ & Peak-hour speed (7--9, 16--19h) & [20, 30] km/h & Est. \\
$v_{\text{mid}}$ & Midday speed (10--15h) & [30, 40] km/h & Est. \\
$v_{\text{off}}$ & Off-peak speed & [35, 50] km/h & Est. \\
\hline
\multicolumn{4}{l}{\textit{Station Parameters}} \\
$\text{Cap}_k$ & Ports per station & [1, 28], $\mu{=}5$ & Data \\
\hline
\multicolumn{4}{l}{\textit{Suitability (Sigmoid)}} \\
$q_c$ & Consumer suitability & [0.003, 0.999], $\mu{=}0.337$ & Data \\
$q_s$ & Provider suitability & [0.004, 0.999], $\mu{=}0.361$ & Data \\
\hline
\multicolumn{4}{l}{\textit{Mobility Cost ($f_3$)}} \\
$\alpha$ & Distance cost weight & 0.7 & Set \\
$\beta$ & Wait time cost weight & 0.3 & Set \\
$W_{\max}$ & Max wait time & 0.5 h & Set \\
\hline
\multicolumn{4}{l}{\scriptsize Data = from dataset; Set = design choice; Est. = estimated; $\mu$ = mean} \\
\hline
\end{tabular}
\end{table}

\subsection{Structure of the Feasible Opportunity Space}
\label{sec:omega_analysis}

Before comparing algorithmic performance, we first examine the structure of the feasible opportunity set $\Omega$. Although the total size of $\Omega$ grows with the number of events, its distribution across EVs is highly uneven. While many EV decision events possess multiple feasible charging sessions across different stations and time slots, a non-negligible subset has only one or two feasible opportunities due to tight waiting tolerances or limited spatial reach.

This asymmetry creates a distinction between flexible EV decision events with many alternatives and fragile EV decision events whose feasibility depends on preserving a small number of critical sessions. As a result, local greedy allocation can inadvertently block fragile EV decision events even when aggregate feasibility exists, motivating the need for globally coordinated and efficiency-aware matching.


\subsection{Comparison with Baseline Methods}

To assess the effectiveness of the proposed NSGA-II, we compare it against three representative baseline methods, namely PPO, Greedy, and Random. These methods capture distinct optimization behaviors: PPO serves as a learning-based baseline, Greedy provides a heuristic baseline based on immediate assignment decisions, and Random offers an unguided reference point. All algorithms are evaluated under the same problem setting with 150 consumer EV decision events and 150 provider EV decision events and are assessed using four objectives: matched provider--consumer decision-event pairs ($f_1$), transferred energy ($f_2$), mobility cost ($f_3$), and suitability score ($f_4$). Because NSGA-II is a multi-objective optimizer that returns a Pareto front rather than a single solution, we first report a balanced solution for direct comparison with the baselines, and then present additional Pareto-front solutions to demonstrate the range of trade-offs captured by the proposed method.



\subsubsection{Comparative Performance of the Proposed and Baseline Algorithms}

\begin{figure}[h]
\centering
\includegraphics[width=0.5\textwidth]{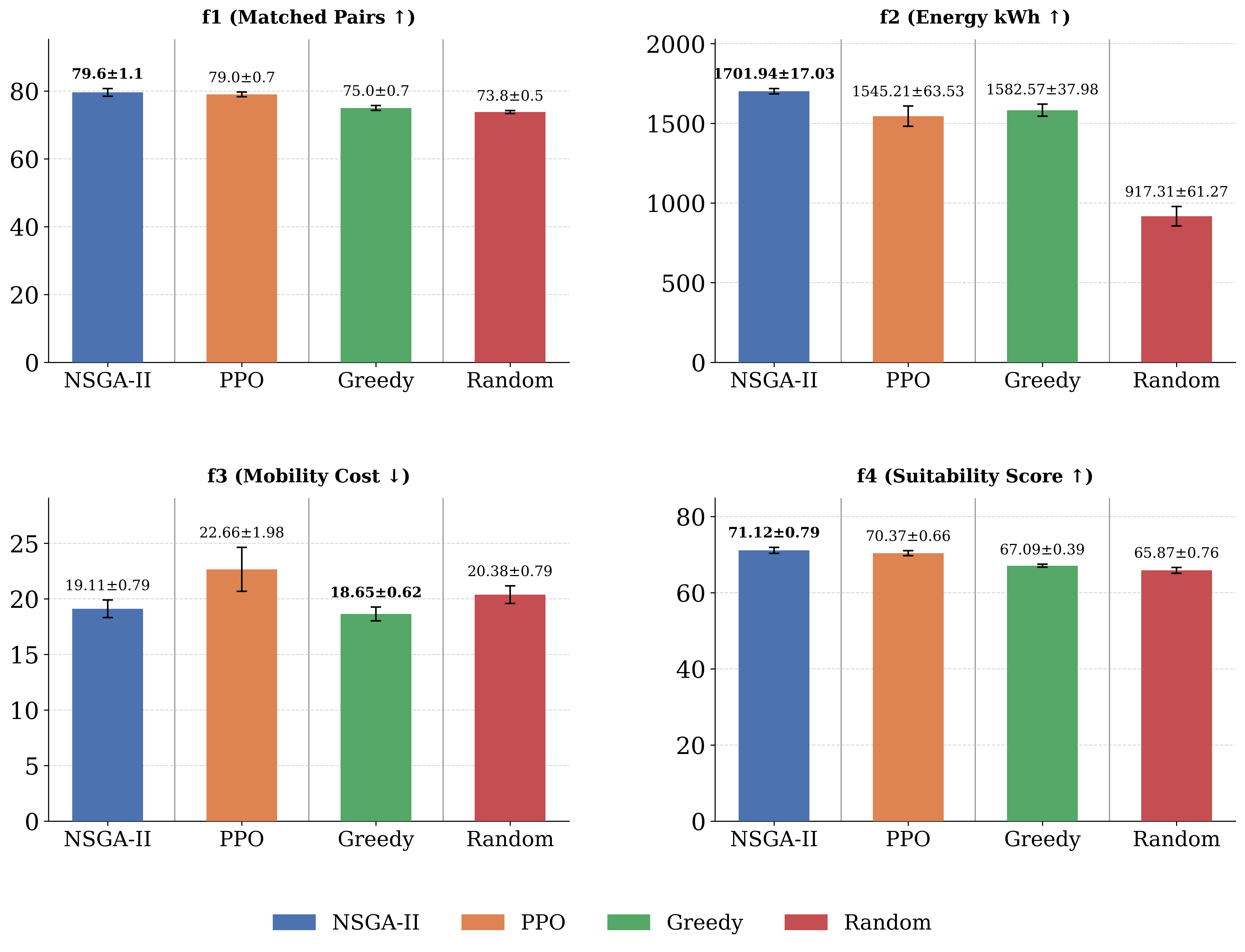}
\caption{Baseline comparison across the four optimization objectives. Bars show the mean objective value over five independent runs, and error bars show the corresponding standard deviation.}
\label{fig:match-baseline-1}
\end{figure}

Table~\ref{tab:main_comparison} and Figure~\ref{fig:match-baseline-1} compare the solutions obtained by NSGA-II, PPO, Greedy, and Random across four objectives. Since the algorithms include stochastic components, the results are reported as \emph{mean $\pm$ standard deviation} over five independent runs.

The results show that PPO and Greedy are strong baselines under this experimental setting, with performance values close to the proposed method on several objectives. 
NSGA-II obtains the highest average number of matched EV--EV trading sessions, with $79.60 \pm 1.14$ successful sessions, corresponding to a matching rate of $53.07 \pm 0.76\%$. This represents a modest improvement over PPO, which achieves $79.00 \pm 0.71$ matches and a matching rate of $52.67 \pm 0.47\%$. The matching improvement over PPO is modest, whereas the transferred-energy and mobility-cost differences are more pronounced. Compared with Greedy and Random, NSGA-II improves the average matching rate by 3.07 and 3.87 percentage points, respectively.

For transferred energy, NSGA-II achieves the highest average value of $1701.94 \pm 17.03$~kWh, compared with $1545.21 \pm 63.53$~kWh for PPO, $1582.57 \pm 37.98$~kWh for Greedy, and $917.31 \pm 61.27$~kWh for Random. This corresponds to improvements of approximately 10.14\%, 7.54\%, and 85.54\%, respectively. Similarly, NSGA-II achieves the highest suitability score of $71.12 \pm 0.79$, improving over PPO, Greedy, and Random by approximately 1.07\%, 6.01\%, and 7.97\%, respectively. These results indicate that NSGA-II provides consistent gains in transferred energy and suitability while also maintaining the highest average matching coverage.

For mobility cost, lower values are preferred. Greedy achieves the lowest average mobility cost, with $18.65 \pm 0.62$, followed by NSGA-II with $19.11 \pm 0.79$, Random with $20.38 \pm 0.79$, and PPO with $22.66 \pm 1.98$. Therefore, NSGA-II does not achieve the best mobility cost, but it remains close to Greedy while substantially improving matching coverage, transferred energy, and suitability. Compared with Greedy, NSGA-II has an approximately 2.47\% higher mobility cost, whereas compared with PPO and Random it reduces mobility cost by approximately 15.67\% and 6.23\%, respectively.

\begin{table*}[t]
\centering
\caption{Overall performance comparison of the proposed method against baseline algorithms. Results are reported as mean $\pm$ standard deviation over five independent runs.}
\label{tab:main_comparison}
\begin{tabular}{lccccc}
\toprule
\textbf{Algorithm} & \textbf{$f_1$ (matches$\uparrow$)} & \textbf{Matching (\%)} & \textbf{$f_2$ (energy kWh$\uparrow$)} & \textbf{$f_3$ (mobility$\downarrow$)} & \textbf{$f_4$ (suitability$\uparrow$)} \\
\midrule
NSGA-II (Proposed) & $79.60 \pm 1.14$ & $53.07 \pm 0.76$ & $1701.94 \pm 17.03$ & $19.11 \pm 0.79$ & $71.12 \pm 0.79$ \\
PPO                & $79.00 \pm 0.71$ & $52.67 \pm 0.47$ & $1545.21 \pm 63.53$ & $22.66 \pm 1.98$ & $70.37 \pm 0.66$ \\
Greedy             & $75.00 \pm 0.71$ & $50.00 \pm 0.47$ & $1582.57 \pm 37.98$ & $18.65 \pm 0.62$ & $67.09 \pm 0.39$ \\
Random             & $73.80 \pm 0.45$ & $49.20 \pm 0.30$ & $917.31 \pm 61.27$ & $20.38 \pm 0.79$ & $65.87 \pm 0.76$ \\
\bottomrule
\end{tabular}
\end{table*}


Taken together, the comparison shows that NSGA-II provides a balanced outcome among the evaluated methods. Its improvement over PPO in matching coverage is modest, but it transfers substantially more energy and achieves lower mobility cost. Greedy remains competitive in mobility cost and energy per match, but it produces fewer matches and lower suitability. 
The main contribution of NSGA-II is its ability to improve several system-level objectives simultaneously while making the mobility-cost trade-off explicit. This trade-off behavior is examined further through the Pareto-front and hypervolume analyses, which show how different objective balances can be selected according to operational priorities.

\subsubsection{Pareto-Front Convergence and Trade-off Analysis of the Proposed and Baseline Algorithms}

\begin{figure*}[h]
\centering
\includegraphics[width=\textwidth]{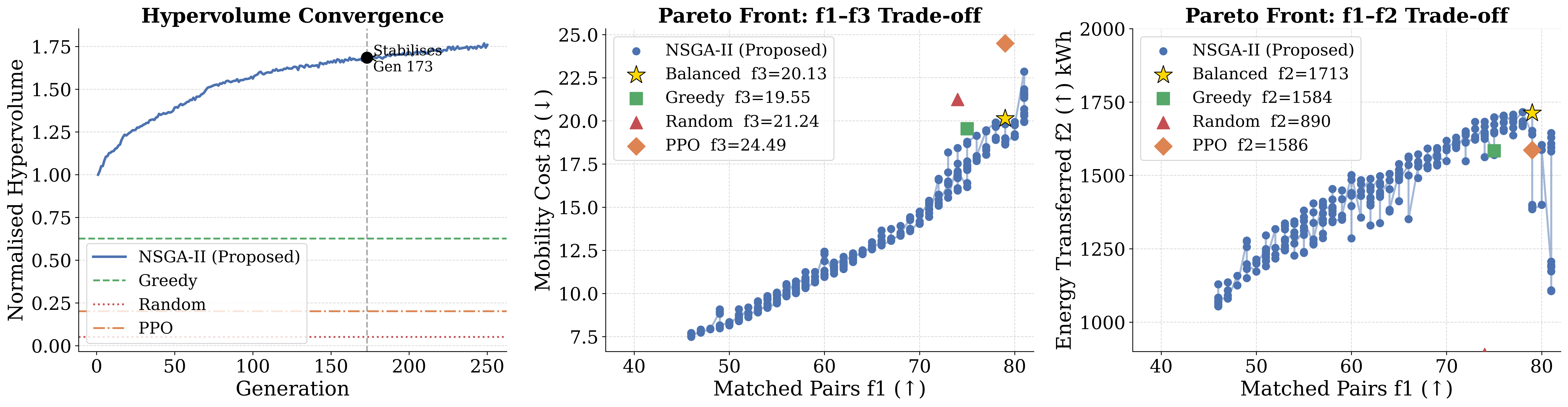}
\caption{Hypervolume convergence and Pareto-front trade-offs of NSGA-II against baseline methods for the representative seed~42 run. The left panel shows the hypervolume progression of NSGA-II across generations. The middle and right panels show the $f_1$--$f_3$ and $f_1$--$f_2$ trade-offs, respectively, with the NSGA-II solution and the baseline methods highlighted for comparison. The hypervolume axis is zero-anchored, while the Pareto-front panels use data-driven axis ranges to preserve the visible shape of the trade-off curve.}
\label{fig:pareto_convergence_tradeoff}
\end{figure*}

Figure~\ref{fig:pareto_convergence_tradeoff} and Table~\ref{tab:pareto_coverage_baseline} provide further insight into the convergence behavior and trade-off characteristics of the proposed NSGA-II compared with the baseline algorithms. We use the hypervolume indicator to assess the convergence and spread of the obtained Pareto front~\cite{zitzler1999multiobjective}. The aggregate baseline comparison is reported as \emph{mean $\pm$ standard} deviation over five independent runs in Table~\ref{tab:main_comparison}. In contrast, the Pareto-front plots and Table~\ref{tab:pareto_coverage_baseline} are reported for the representative seed~42 run, because each NSGA-II row corresponds to a distinct feasible solution selected from the final Pareto front and should therefore not be averaged.

The left panel of Figure~\ref{fig:pareto_convergence_tradeoff} shows the hypervolume progression of NSGA-II over generations. The hypervolume increases steadily during the search process and stabilizes around Generation~173, indicating that the population converges toward a stable approximation of the Pareto front. The baseline methods are shown as fixed hypervolume levels because each baseline produces a single matching solution rather than an evolving Pareto population. These fixed levels remain below the final hypervolume achieved by NSGA-II, indicating that the proposed method explores a broader and higher-quality trade-off region than PPO, Greedy, and Random.

The middle and right panels of Figure~\ref{fig:pareto_convergence_tradeoff} illustrate the trade-offs captured by the Pareto front in the $f_1$--$f_3$ and $f_1$--$f_2$ objective spaces, respectively. In the $f_1$--$f_3$ projection, lower mobility cost is obtained at the expense of fewer matched pairs, while higher matching performance generally requires higher mobility cost. In the $f_1$--$f_2$ projection, increases in matched pairs are generally associated with higher transferred energy, although the front also reveals local variability due to the interaction among multiple objectives and feasibility constraints. These analyses confirm that the optimization problem is inherently multi-objective, with no single solution simultaneously optimizing all objectives.


The balanced NSGA-II solution highlighted in Figure~\ref{fig:pareto_convergence_tradeoff} represents a compromise region of the Pareto front, achieving high matching coverage and transferred energy while avoiding an excessive increase in mobility cost. In the representative seed~42 run, this solution achieves 78 matched pairs, corresponding to a matching rate of 52.00\%, with approximately 1713~kWh of transferred energy and a mobility cost of approximately 20.13. Compared with PPO, Greedy, and Random, the balanced NSGA-II solution provides a competitive matching level while improving the trade-off between transferred energy and mobility cost. Although Greedy obtains a slightly lower mobility cost, it does so with fewer matched pairs and lower suitability. Similarly, PPO achieves a comparable number of matches, but with substantially higher mobility cost. This behavior is consistent with the role of NSGA-II as a multi-objective optimizer that prioritizes balanced compromise solutions rather than single-objective extremal points.

\begin{table*}[t]
\centering
\caption{Pareto-front coverage of NSGA-II across different objectives for the representative seed~42 run. Each NSGA-II row represents a distinct solution selected from the Pareto front.}
\label{tab:pareto_coverage_baseline}
\begin{tabular}{lccccc}
\toprule
\textbf{Algorithm} & \textbf{$f_1$ (matches$\uparrow$)} & \textbf{Matching (\%)} & \textbf{$f_2$ (energy kWh$\uparrow$)} & \textbf{$f_3$ (mobility$\downarrow$)} & \textbf{$f_4$ (suitability$\uparrow$)} \\
\midrule
NSGA-II best-$f_1$ ($\uparrow$) & 81 & 54.00 & 1644.86 & 20.6810 & 72.4823 \\
NSGA-II best-$f_2$ ($\uparrow$) & 78 & 52.00 & 1716.04 & 19.9124 & 70.0830 \\
NSGA-II best-$f_3$ ($\downarrow$) & 46 & 30.67 & 1057.84 & 7.4857 & 42.4911 \\
NSGA-II best-$f_4$ ($\uparrow$) & 81 & 54.00 & 1105.89 & 21.4806 & 72.7895 \\
PPO                              & 79 & 52.67 & 1586.45 & 24.4932 & 70.6471 \\
Greedy                           & 75 & 50.00 & 1583.98 & 19.5476 & 67.1144 \\
Random                           & 74 & 49.33 & 890.19 & 21.2445 & 65.5223 \\
\bottomrule
\end{tabular}
\end{table*}

Table~\ref{tab:pareto_coverage_baseline} further demonstrates the coverage of the Pareto front by reporting representative NSGA-II solutions optimized toward different objectives. The best-$f_1$ solution achieves 81 matches, corresponding to a matching rate of 54.00\%, which is higher than all baseline methods in the representative run. This solution also achieves high transferred energy of 1644.86~kWh and the highest suitability among the high-matching solutions, but with a higher mobility cost than Greedy. The best-$f_2$ solution achieves the highest transferred energy, 1716.04~kWh, while maintaining 78 matches and a mobility cost of 19.9124, which is close to Greedy. The best-$f_3$ solution minimizes mobility cost to 7.4857, but this comes with a substantial reduction in matched pairs, transferred energy, and suitability, clearly illustrating the cost of optimizing mobility in isolation. The best-$f_4$ solution attains the highest suitability score of 72.7895 together with 81 matches, but transfers only 1105.89~kWh, showing that suitability-oriented matching may reduce transferred energy.

Compared with the baseline methods, the NSGA-II Pareto front provides a wider range of feasible operating points. PPO achieves 79 matches and strong suitability, but has the highest mobility cost among the baseline methods. Greedy provides the lowest baseline mobility cost and competitive energy transfer, but produces fewer matches and lower suitability than the high-coverage NSGA-II solutions. Random performs reasonably in matching coverage for this instance, but transfers substantially less energy and has lower suitability. These results show that NSGA-II is not limited to a single fixed trade-off; instead, it provides decision-makers with multiple feasible solutions that can be selected according to operational priorities.

Taken together, the results confirm two key points. First, the proposed NSGA-II converges reliably toward a stable Pareto front with greater hypervolume than the baseline methods. Second, it offers a richer and more effective set of trade-off solutions, including a balanced solution that outperforms the baseline algorithms on the most important system-level objectives while preserving acceptable mobility cost. This supports the suitability of the proposed NSGA-II for EV--EV energy trading scenarios where multiple conflicting objectives must be considered jointly.

\subsection{Comparison with State-of-the-art Methods}

To assess the effectiveness of the proposed NSGA-II against existing approaches, we compare it with two representative state-of-the-art methods, namely BMNN and DAMA. These methods reflect different optimization paradigms reported in existing EV--EV energy trading studies. BMNN serves as a recent optimization-based matching approach, while DAMA represents an auction-driven coordination mechanism for provider--consumer matching. All algorithms are evaluated under the same problem setting with 150 consumer EV decision events and 150 provider EV decision events and are assessed using four objectives: matched pairs ($f_1$), transferred energy ($f_2$), mobility cost ($f_3$), and suitability score ($f_4$). Since NSGA-II is a multi-objective optimizer that produces a Pareto front rather than a single solution, the balanced NSGA-II solution from each run is used for direct comparison with BMNN and DAMA.

\subsubsection{Performance of the Proposed NSGA-II Against State-of-the-Art Methods}


\begin{figure}[h]
\centering
\includegraphics[width=0.5\textwidth]{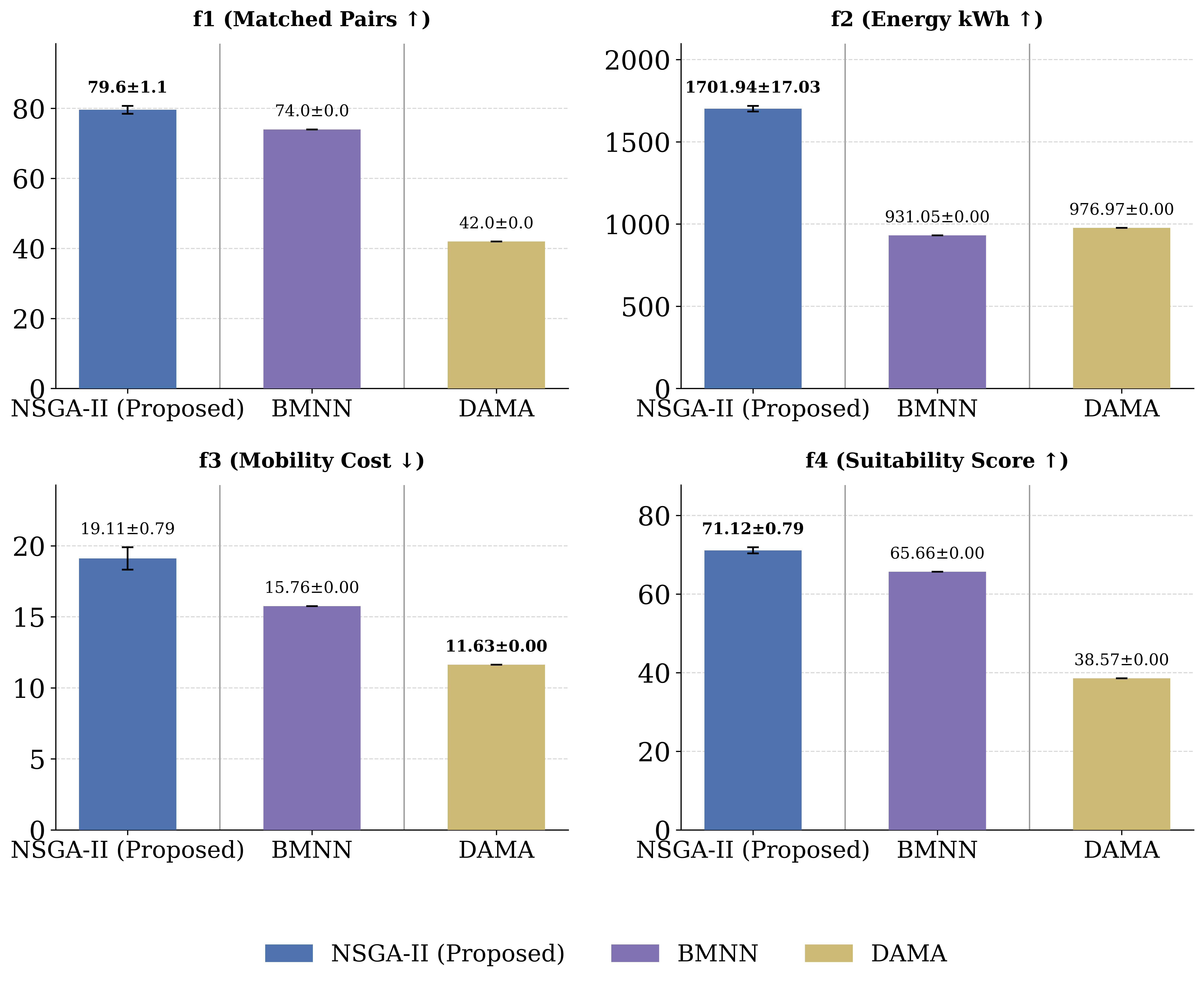}
\caption{Performance comparison of NSGA-II against state-of-the-art EV--EV matching methods. Bars show the mean objective value over five runs, and error bars show the corresponding standard deviation.}
\label{fig:match-sota}
\end{figure}

\begin{table*}[t]
\centering
\caption{Comparative performance of state-of-the-art methods across matching, energy, mobility, and suitability objectives. Results are reported as mean $\pm$ standard deviation over five runs. BMNN and DAMA are deterministic for the fixed experimental instance; therefore, their repeated-run standard deviations are zero.}
\label{tab:total_objectives_bmnn}
\begin{tabular}{lccccc}
\toprule
\textbf{Algorithm} & \textbf{$f_1$ (matches$\uparrow$)} & \textbf{Matching (\%)} & \textbf{$f_2$ (energy kWh$\uparrow$)} & \textbf{$f_3$ (mobility$\downarrow$)} & \textbf{$f_4$ (suitability$\uparrow$)} \\
\midrule
NSGA-II (Proposed) & $79.60 \pm 1.14$ & $53.07 \pm 0.76$ & $1701.94 \pm 17.03$ & $19.11 \pm 0.79$ & $71.12 \pm 0.79$ \\
BMNN               & $74.00 \pm 0.00$ & $49.33 \pm 0.00$ & $931.05 \pm 0.00$  & $15.76 \pm 0.00$ & $65.66 \pm 0.00$ \\
DAMA               & $42.00 \pm 0.00$ & $28.00 \pm 0.00$ & $976.97 \pm 0.00$  & $11.63 \pm 0.00$ & $38.57 \pm 0.00$ \\
\bottomrule
\end{tabular}
\end{table*}

Figure~\ref{fig:match-sota} and Table~\ref{tab:total_objectives_bmnn} compare the proposed NSGA-II with BMNN and DAMA across four objectives. Results are reported as \emph{mean $\pm$ standard deviation} over five runs. The proposed NSGA-II is stochastic and is therefore evaluated over five independent random seeds. BMNN and DAMA are deterministic for the fixed input instance; hence, repeated executions produce identical outputs and their standard deviations are zero.

The results show that NSGA-II achieves stronger system-level outcomes than BMNN and DAMA in terms of matching coverage, transferred energy, and charging-node suitability. NSGA-II obtains $79.60 \pm 1.14$ successful matches, corresponding to a matching rate of $53.07 \pm 0.76\%$. In comparison, BMNN achieves $74.00 \pm 0.00$ matches with a matching rate of $49.33 \pm 0.00\%$, while DAMA achieves $42.00 \pm 0.00$ matches with a matching rate of $28.00 \pm 0.00\%$. Thus, NSGA-II improves the average matching rate by 3.74 percentage points over BMNN and 25.07 percentage points over DAMA.

For transferred energy, NSGA-II achieves the highest average value of $1701.94 \pm 17.03$~kWh, compared with $931.05 \pm 0.00$~kWh for BMNN and $976.97 \pm 0.00$~kWh for DAMA. This corresponds to improvements of approximately 82.80\% over BMNN and 74.21\% over DAMA. NSGA-II also achieves the highest suitability score, with $71.12 \pm 0.79$, compared with $65.66 \pm 0.00$ for BMNN and $38.57 \pm 0.00$ for DAMA. These results indicate that the proposed method produces more suitable provider--consumer and charging-node assignments while transferring substantially more energy.

For mobility cost, lower values are preferred. DAMA achieves the lowest mobility cost, with $11.63 \pm 0.00$, followed by BMNN with $15.76 \pm 0.00$, while NSGA-II records $19.11 \pm 0.79$. This result shows that the proposed method accepts higher travel-related cost in exchange for stronger matching coverage, transferred energy, and charging-node suitability. This behavior is consistent with the role of NSGA-II as a multi-objective optimizer that seeks balanced system-level trade-offs rather than minimizing travel overhead alone.

It is also important to clarify that the lower performance of BMNN and DAMA on matching coverage, transferred energy, and suitability does not imply that these methods are generally inferior in all EV--EV trading settings. Rather, it reflects differences in optimization alignment under the studied problem setting. BMNN prioritizes proximity-oriented matching, while DAMA focuses on auction-based clearing. In contrast, the proposed NSGA-II explicitly optimizes the four objectives under the same journey-aware feasibility constraints. As a result, BMNN and DAMA retain an advantage in mobility cost, but they achieve lower matching coverage, transferred energy, and suitability under the objective structure considered in this paper.

In summary, the comparison shows that NSGA-II is more effective for the prediction-guided, journey-aware multi-objective matching problem considered in this study. Relative to BMNN and DAMA, it improves matching coverage, transferred energy, and charging-node suitability, although at the expense of higher mobility cost. These results support the value of NSGA-II in EV--EV energy trading scenarios where the goal is to balance multiple competing objectives rather than optimize only proximity, auction surplus, or travel efficiency.

\subsubsection{Pareto-Front Convergence and Trade-off Analysis Against State-of-the-Art Methods}

\begin{figure*}[h]
\centering
\includegraphics[width=\textwidth]{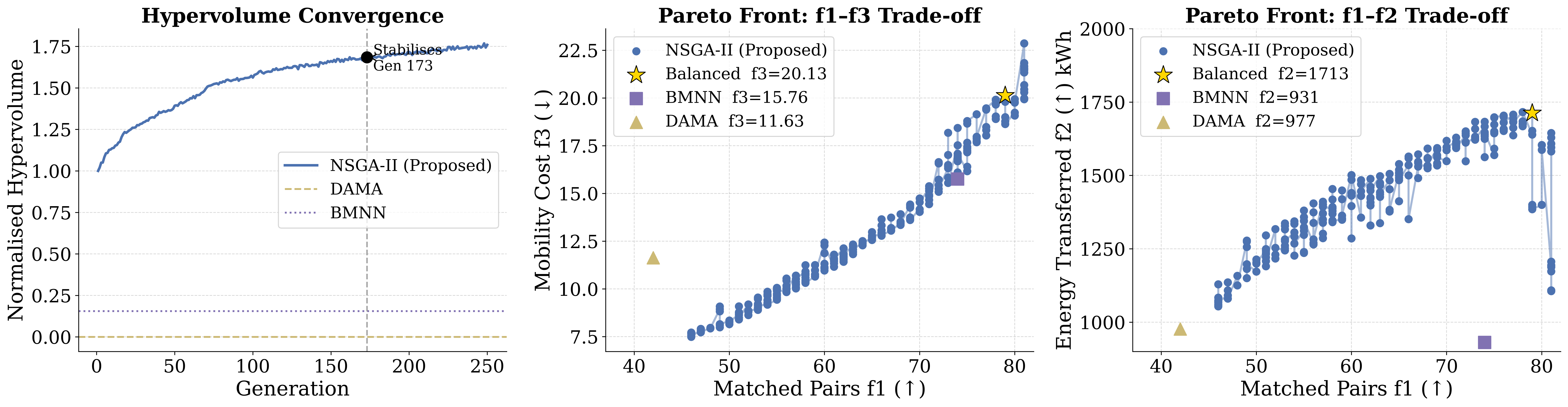}
\caption{Hypervolume convergence and Pareto-front trade-offs of NSGA-II against state-of-the-art methods for the representative seed~42 run. The left panel shows the hypervolume progression of NSGA-II across generations. The middle and right panels show the $f_1$--$f_3$ and $f_1$--$f_2$ trade-offs, respectively, with the NSGA-II solution and the state-of-the-art reference methods highlighted for comparison. The hypervolume axis is zero-anchored, while the Pareto-front panels use data-driven axis ranges to preserve the visible shape of the trade-off curve.}
\label{fig:match-sota-convergence}
\end{figure*}

Figure~\ref{fig:match-sota-convergence} and Table~\ref{tab:pareto_coverage_bmnn_dama_nsga} provide a deeper comparison between the proposed NSGA-II and the state-of-the-art methods, BMNN and DAMA, by examining both convergence behavior and Pareto-front trade-offs. The aggregate performance comparison in Table~\ref{tab:total_objectives_bmnn} is reported as $mean \pm standard deviation$ over five runs. In contrast, the Pareto-front plots and Table~\ref{tab:pareto_coverage_bmnn_dama_nsga} are reported for the representative seed~42 run, because each NSGA-II row corresponds to a distinct feasible solution selected from the final Pareto front and should therefore not be averaged.

The left panel of Figure~\ref{fig:match-sota-convergence} shows the hypervolume progression of NSGA-II over generations. The hypervolume increases steadily during the search process and stabilizes around Generation~173, indicating that the population converges toward a stable approximation of the Pareto front. BMNN and DAMA are shown as fixed hypervolume levels because each method produces a single matching solution rather than an evolving Pareto population. These fixed levels remain below the final hypervolume achieved by NSGA-II, indicating that the proposed method explores a broader and higher-quality trade-off region than the state-of-the-art reference methods.

\begin{table*}[t]
\centering
\caption{Comparison of representative Pareto-front solutions of NSGA-II across different objectives, with BMNN and DAMA as references. Each NSGA-II row represents a distinct solution selected from the representative seed~42 Pareto front.}
\label{tab:pareto_coverage_bmnn_dama_nsga}
\begin{tabular}{lccccc}
\toprule
\textbf{Algorithm} & \textbf{$f_1$ (matches$\uparrow$)} & \textbf{Matching (\%)} & \textbf{$f_2$ (energy kWh$\uparrow$)} & \textbf{$f_3$ (mobility$\downarrow$)} & \textbf{$f_4$ (suitability$\uparrow$)} \\
\midrule
NSGA-II best-$f_1$ ($\uparrow$) & 81 & 54.00 & 1644.86 & 20.6810 & 72.4823 \\
NSGA-II best-$f_2$ ($\uparrow$) & 78 & 52.00 & 1716.04 & 19.9124 & 70.0830 \\
NSGA-II best-$f_3$ ($\downarrow$) & 46 & 30.67 & 1057.84 & 7.4857 & 42.4911 \\
NSGA-II best-$f_4$ ($\uparrow$) & 81 & 54.00 & 1105.89 & 21.4806 & 72.7895 \\
BMNN                              & 74 & 49.33 & 931.05 & 15.7566 & 65.6634 \\
DAMA                              & 42 & 28.00 & 976.97 & 11.6294 & 38.5708 \\
\bottomrule
\end{tabular}
\end{table*}

The middle and right panels of Figure~\ref{fig:match-sota-convergence} illustrate the trade-offs captured by the Pareto front in the $f_1$--$f_3$ and $f_1$--$f_2$ objective spaces, respectively. In the $f_1$--$f_3$ projection, lower mobility cost is obtained at the expense of fewer matched pairs, while higher matching performance generally requires higher mobility cost. In the $f_1$--$f_2$ projection, higher numbers of matched pairs are generally associated with higher transferred energy, although local variability remains because all four objectives are optimized simultaneously under feasibility constraints. These analyses confirm that the EV--EV energy trading problem is inherently multi-objective, with no single solution optimizing all objectives at once.

The balanced NSGA-II solution highlighted in Figure~\ref{fig:match-sota-convergence} lies in a favorable trade-off region of the Pareto front. In the representative seed~42 run, it achieves 78 matched pairs, corresponding to a matching rate of 52.00\%, with approximately 1713~kWh of transferred energy and a mobility cost of approximately 20.13. Compared with BMNN and DAMA, this solution provides stronger matching coverage and substantially higher transferred energy, while maintaining a reasonable mobility-cost trade-off. BMNN and DAMA achieve lower mobility costs, but they do so with fewer matches, lower transferred energy, and lower suitability. This behavior is consistent with the design goal of NSGA-II, which is to identify balanced compromise solutions across conflicting objectives rather than optimize a single objective in isolation.

Table~\ref{tab:pareto_coverage_bmnn_dama_nsga} further demonstrates the breadth of the Pareto front by reporting representative NSGA-II solutions optimized toward different objectives. The best-$f_1$ solution achieves 81 matches, corresponding to a matching rate of 54.00\%, which is higher than both BMNN and DAMA. The best-$f_2$ solution achieves the highest transferred energy, 1716.04~kWh, while maintaining 78 matches and a mobility cost of 19.9124. The best-$f_3$ solution reduces mobility cost to 7.4857, but this comes with a substantial reduction in matched pairs, transferred energy, and suitability, clearly illustrating the loss incurred when mobility cost is optimized too aggressively. The best-$f_4$ solution achieves the highest suitability score of 72.7895 with 81 matches, although with lower transferred energy than the best-$f_1$ and best-$f_2$ solutions. These results show that the Pareto front produced by NSGA-II spans a wide range of feasible trade-offs and gives decision-makers flexibility to select solutions according to system priorities.

Compared with BMNN and DAMA, the NSGA-II Pareto front provides stronger coverage of high-matching, high-energy, and high-suitability operating points. BMNN achieves 74 matches and a mobility cost of 15.7566, but transfers only 931.05~kWh and achieves a suitability score of 65.6634. DAMA achieves the lowest mobility cost among the reference methods, 11.6294, but only obtains 42 matches and a suitability score of 38.5708. Therefore, BMNN and DAMA remain effective when the priority is mobility-cost reduction, whereas NSGA-II provides a broader set of trade-off solutions for balancing matching coverage, transferred energy, mobility cost, and suitability jointly.

Taken together, the results confirm two key findings. First, NSGA-II converges reliably toward a stable Pareto front with substantially greater hypervolume than BMNN and DAMA. Second, it offers a richer and more effective set of trade-off solutions, including high-coverage, high-energy, low-mobility, and high-suitability alternatives. This demonstrates the advantage of NSGA-II for EV--EV energy trading scenarios where multiple conflicting objectives must be considered jointly.

\subsection{Ablation Study of Prediction-Derived Suitability Guidance}

Table~\ref{tab:ablation_suitability_guidance} evaluates the contribution of prediction-derived suitability guidance in the proposed NSGA-II. In the ablated variant, prediction-derived suitability is removed from the objective vector, initialization, mutation, repair, and balanced-solution selection. The suitability score $f_4$ is computed only after optimization as an evaluation metric.

The results show that removing suitability guidance does not substantially reduce matching coverage, which decreases only from $53.07\%$ to $52.67\%$. This indicates that the feasible opportunity set $\Omega$ still contains sufficient spatially and temporally feasible provider--consumer sessions. However, transferred energy decreases substantially, from $1701.94$~kWh to $1407.70$~kWh, and the average transferred energy per match decreases from $21.39$ to $17.82$~kWh. This suggests that prediction-derived suitability guidance contributes mainly to selecting higher-quality trading sessions, rather than merely increasing the number of feasible matches.

\begin{table*}[t]
\centering
\caption{Ablation study of prediction-derived suitability guidance in the proposed NSGA-II. NSGA-II-NoSuitabilityGuidance removes prediction-derived suitability from the objective vector, initialization, mutation, repair, and balanced-solution selection. The suitability score $f_4$ is reported post hoc for the ablated variant. Results are reported as mean $\pm$ standard deviation over five seeds.}
\label{tab:ablation_suitability_guidance}
\begin{tabular}{lccccc}
\toprule
\textbf{Algorithm} & \textbf{$f_1$ (matches$\uparrow$)} & \textbf{Matching (\%)} & \textbf{$f_2$ (energy kWh$\uparrow$)} & \textbf{$f_3$ (mobility$\downarrow$)} & \textbf{$f_4$ (suitability$\uparrow$)} \\
\midrule
NSGA-II (Proposed)                  & $79.60 \pm 1.14$ & $53.07 \pm 0.76$ & $1701.94 \pm 17.03$ & $19.11 \pm 0.79$ & $71.12 \pm 0.79$ \\
NSGA-II-NoSuitabilityGuidance   & $79.00 \pm 1.41$ & $52.67 \pm 0.94$ & $1407.70 \pm 34.92$ & $18.96 \pm 1.24$ & $70.32 \pm 1.38$ \\
\bottomrule
\end{tabular}
\end{table*}

The post hoc suitability score also decreases from $71.12$ to $70.32$, showing a modest reduction in prediction alignment. The relatively small change in $f_4$ is expected because all variants operate on the same feasibility-filtered opportunity set $\Omega$. Nevertheless, the substantial reduction in transferred energy confirms that suitability guidance helps steer the search toward more effective provider--consumer--charging-node assignments within the feasible opportunity space.


\subsection{Computational Efficiency and Runtime Scalability}

\begin{figure}[h]
\centering
\includegraphics[width=0.5\textwidth]{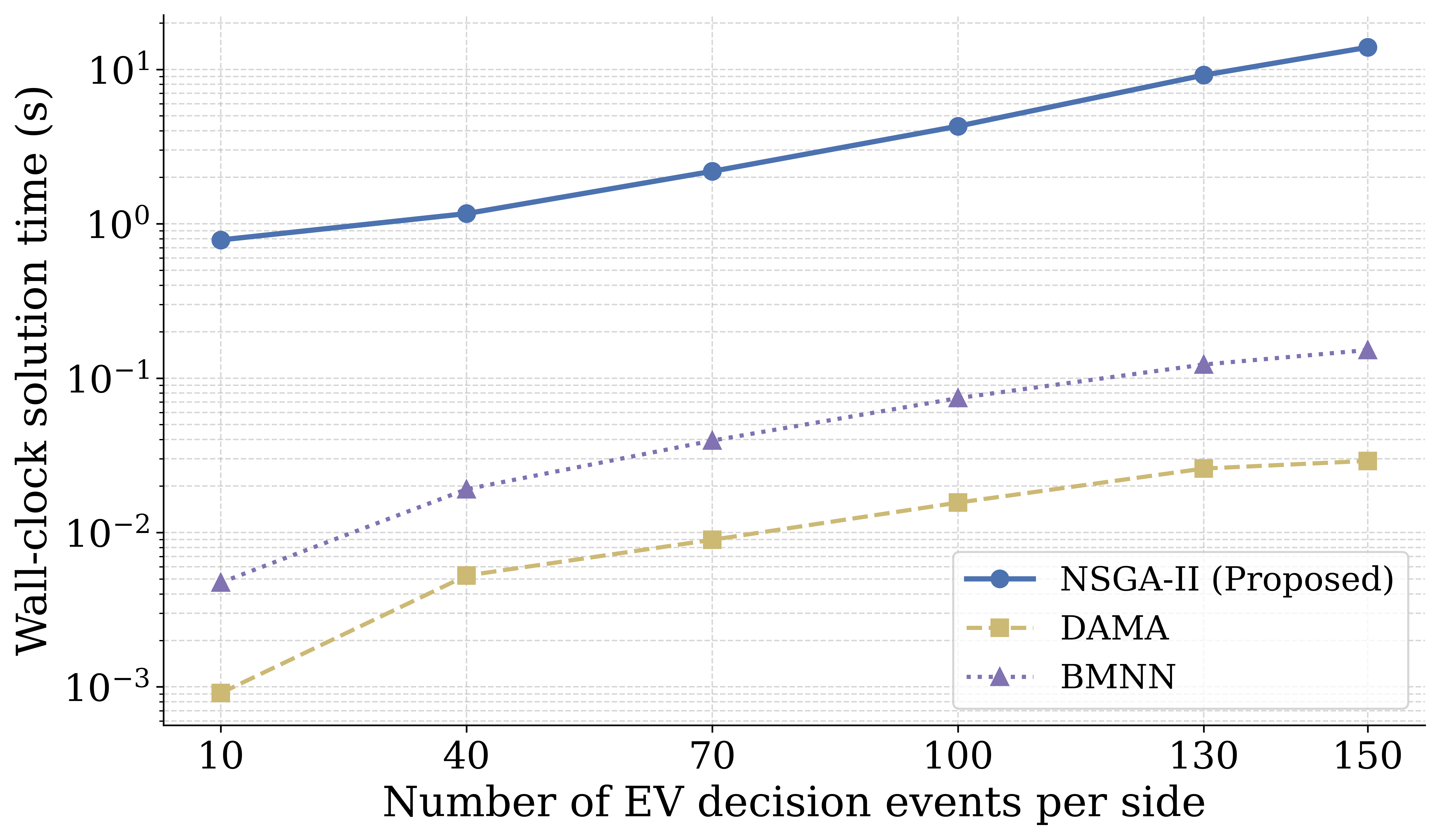}
\caption{Runtime comparison of NSGA-II, DAMA, and BMNN as the number of EV decision events per side increases. Each setting contains an equal number of consumer and provider EV decision events; for example, 150 corresponds to 150 consumer EV decision events and 150 provider EV decision events. The y-axis is shown on a logarithmic scale.}
\label{fig:time}
\end{figure}

Figure~\ref{fig:time} compares the solution times of NSGA-II, DAMA, and BMNN as the number of EV decision events per side increases. This runtime scalability analysis is reported for the representative seed~42 run, whereas the main performance comparison reports $mean \pm standard deviation$ over five independent seeds. The log-scale y-axis shows that all methods require more computation as the problem size grows, but the proposed NSGA-II incurs substantially higher runtime than the two baselines. Specifically, NSGA-II increases from approximately 0.7~s for 10 EV decision events per side to about 20~s for 150 EV decision events per side, whereas DAMA and BMNN remain in the millisecond range across all settings. This behavior is expected because NSGA-II performs iterative multi-objective evolutionary search to explore trade-offs among competing objectives, while DAMA and BMNN rely on less computationally intensive decision procedures. Therefore, the proposed method introduces a clear runtime overhead in exchange for higher-quality optimization. From an operational perspective, this overhead remains acceptable provided that the matching process is executed one charging cycle before the actual charging event, allowing sufficient time to compute a feasible and high-quality solution.

\subsection{Slot-wise Matching Performance of NSGA-II}

\begin{figure}[h]
\centering
\includegraphics[width=0.5\textwidth]{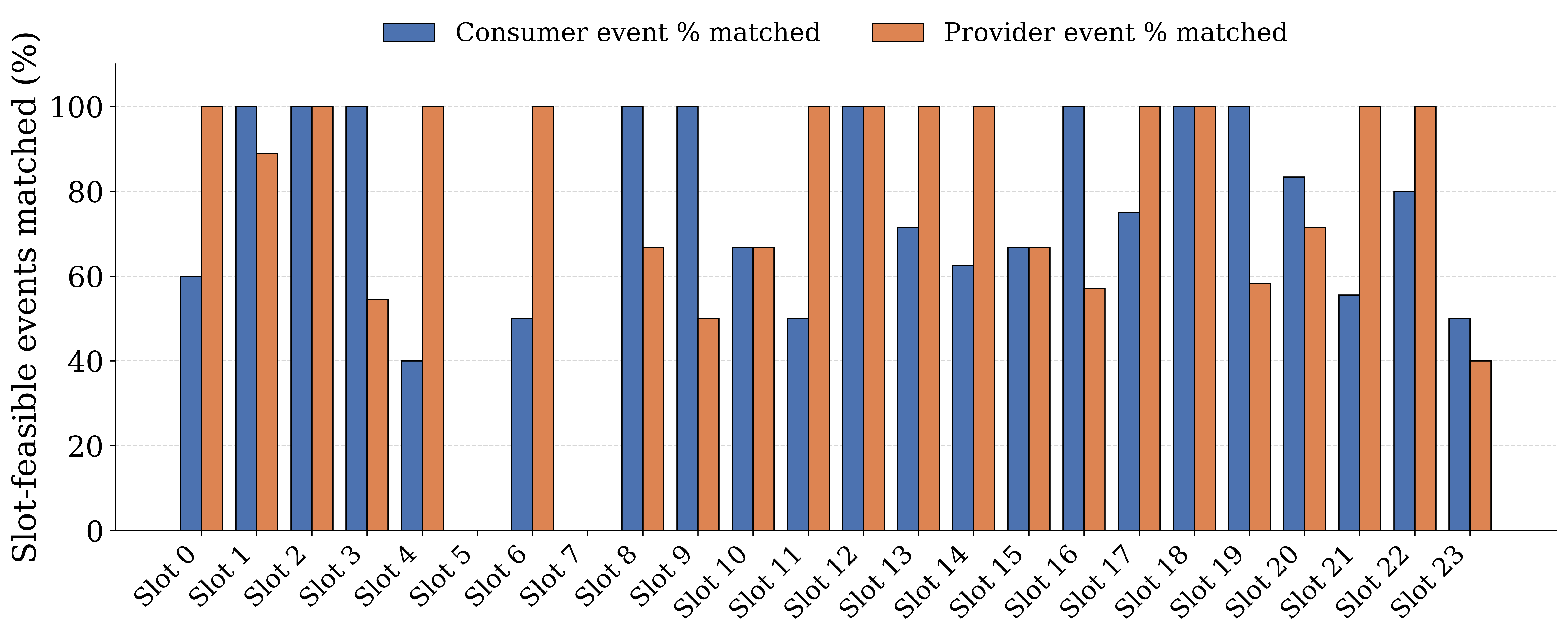}
\caption{Slot-wise feasible-opportunity coverage under NSGA-II for the representative seed~42 run. Consumer and provider percentages are computed relative to the slot-feasible EV decision events that appear in at least one feasible session $(i,j,k,t)\in\Omega$ in the corresponding slot.}
\label{fig:slotwise_matching}
\end{figure}

Table~\ref{tab:slotwise_pg_nsga2} and Figure~\ref{fig:slotwise_matching} present the slot-wise feasible-opportunity coverage of NSGA-II for the representative seed~42 run. The table reports, for each time slot, the number of slot-feasible consumer and provider EV decision events, the corresponding matched decision-event counts, and the conditional matching percentages. A decision event is considered slot-feasible in slot $t$ if it appears in at least one feasible session tuple $(i,j,k,t)\in\Omega$ for that slot. Therefore, the slot-wise denominators vary across time and are generally lower than the 150 sampled consumer and provider EV decision events used in the full optimization instance.

The slot-wise percentages should be interpreted as conditional matching rates within individual time slots, rather than as the aggregate matching coverage reported in the comparative evaluation. The aggregate matching coverage is computed as the total number of selected EV--EV sessions divided by the number of sampled consumer EV decision events in the optimization instance. In the representative seed~42 run, the matched counts sum to 79 selected sessions, corresponding to an aggregate matching coverage of $79/150=52.67\%$. Since each selected session contains one matched consumer EV decision event and one matched provider EV decision event, the matched consumer and provider counts are identical in each slot. By contrast, the slot-wise consumer and provider percentages are computed using slot-specific feasible-event denominators. Thus, the slot-wise results explain the temporal distribution of successful matches, while the aggregate coverage reports the full-instance matching performance.

\begin{table*}[t]
\centering
\caption{Slot-wise feasible and matched consumer--provider EV decision events under NSGA-II for the representative seed~42 run. Slot-wise percentages are conditional on the slot-feasible decision events in each time slot and should not be averaged to obtain the aggregate matching coverage.}
\label{tab:slotwise_pg_nsga2}
\scriptsize
\begin{tabular}{lrrrrrr}
\toprule
\textbf{Slot} &
\textbf{Feasible Consumer} &
\textbf{Matched Consumer} &
\textbf{Consumer Match (\%)} &
\textbf{Feasible Provider} &
\textbf{Matched Provider} &
\textbf{Provider Match (\%)} \\
\midrule
\textbf{Slot 0}  & 5 & 3 & 60.00 & 3  & 3 & 100.00 \\
\textbf{Slot 1}  & 8 & 8 & 100.00 & 9  & 8 & 88.89 \\
\textbf{Slot 2}  & 3 & 3 & 100.00 & 3  & 3 & 100.00 \\
\textbf{Slot 3}  & 6 & 6 & 100.00 & 11 & 6 & 54.55 \\
\textbf{Slot 4}  & 5 & 2 & 40.00 & 2  & 2 & 100.00 \\
\textbf{Slot 5}  & 0 & 0 & 0.00  & 0  & 0 & 0.00 \\
\textbf{Slot 6}  & 2 & 1 & 50.00 & 1  & 1 & 100.00 \\
\textbf{Slot 7}  & 0 & 0 & 0.00  & 0  & 0 & 0.00 \\
\textbf{Slot 8}  & 2 & 2 & 100.00 & 3  & 2 & 66.67 \\
\textbf{Slot 9}  & 1 & 1 & 100.00 & 2  & 1 & 50.00 \\
\textbf{Slot 10} & 3 & 2 & 66.67 & 3  & 2 & 66.67 \\
\textbf{Slot 11} & 4 & 2 & 50.00 & 2  & 2 & 100.00 \\
\textbf{Slot 12} & 1 & 1 & 100.00 & 1  & 1 & 100.00 \\
\textbf{Slot 13} & 7 & 5 & 71.43 & 5  & 5 & 100.00 \\
\textbf{Slot 14} & 8 & 5 & 62.50 & 5  & 5 & 100.00 \\
\textbf{Slot 15} & 3 & 2 & 66.67 & 3  & 2 & 66.67 \\
\textbf{Slot 16} & 4 & 4 & 100.00 & 7  & 4 & 57.14 \\
\textbf{Slot 17} & 8 & 6 & 75.00 & 6  & 6 & 100.00 \\
\textbf{Slot 18} & 3 & 3 & 100.00 & 3  & 3 & 100.00 \\
\textbf{Slot 19} & 7 & 7 & 100.00 & 12 & 7 & 58.33 \\
\textbf{Slot 20} & 6 & 5 & 83.33 & 7  & 5 & 71.43 \\
\textbf{Slot 21} & 9 & 5 & 55.56 & 5  & 5 & 100.00 \\
\textbf{Slot 22} & 5 & 4 & 80.00 & 4  & 4 & 100.00 \\
\textbf{Slot 23} & 4 & 2 & 50.00 & 5  & 2 & 40.00 \\
\bottomrule
\end{tabular}
\end{table*}

The results show substantial temporal variation across the 24 slots, indicating that matching feasibility depends strongly on slot-specific consumer--provider availability, WCP accessibility, and feasibility constraints. Several slots achieve high conditional coverage because most slot-feasible decision events are successfully assigned. On the consumer side, full conditional coverage is observed in Slots~1, 2, 3, 8, 9, 12, 16, 18, and~19. Among these, Slots~1, 3, and~19 are particularly informative because they combine high percentages with comparatively larger feasible-event counts. On the provider side, full conditional coverage is observed in Slots~0, 2, 4, 6, 11, 12, 13, 14, 17, 18, 21, and~22. Provider-side coverage in Slots~13, 14, 17, and~21 is notable because these slots include multiple matched sessions while still fully utilizing the slot-feasible provider set.

Slots~5 and~7 contain no feasible session opportunities in $\Omega$ and therefore produce no matches. Other slots exhibit partial coverage despite having feasible events. For example, Slot~4 matches only 2 out of 5 feasible consumer decision events, while Slot~23 matches 2 out of 5 feasible provider decision events. These cases show that the existence of slot-feasible events does not guarantee assignment, because the final selected sessions must jointly satisfy consumer--provider compatibility, WCP--slot capacity, and the multi-objective trade-offs among coverage, transferred energy, mobility cost, and suitability.

The comparison between consumer and provider percentages also reveals slot-specific demand--supply imbalance. For example, in Slot~19, all 7 feasible consumer decision events are matched, whereas 7 out of 12 feasible provider decision events are used. This suggests a relative provider surplus in that slot. Conversely, Slot~21 matches 5 out of 9 feasible consumer decision events but uses all 5 feasible provider decision events, indicating tighter provider-side availability. Similar asymmetries across slots confirm that matching performance is influenced not only by the number of feasible participants, but also by the balance and compatibility between the two sides.

Taken together, the slot-wise results show that NSGA-II identifies stronger matching opportunities in favourable time slots while maintaining feasibility under varying temporal conditions. The variation across slots highlights the difficulty of EV--EV energy trading in dynamic settings, where feasible opportunities, counterpart availability, and objective trade-offs change over time. These findings provide a fine-grained view of the proposed method and show that the aggregate performance is driven by its ability to exploit high-opportunity slots while handling constrained slots conservatively.

\subsection{Sensitivity Analysis}

Table~\ref{tab:sensitivity} reports the sensitivity analysis of NSGA-II on the 150-consumer/150-provider instance. Each configuration was evaluated over three independent tuning runs, while all other parameters were fixed at their default values during the corresponding sweep. Since the results are averaged over multiple runs, all objective values in the table are reported as average values. The selected configuration is not necessarily the best setting for every individual objective. This is expected because the four objectives are conflicting: settings that improve matching coverage do not always improve transferred energy, mobility cost, or suitability. Therefore, the final hyperparameter setting was selected using a balanced multi-objective criterion rather than a single-objective maximum. In particular, the selected settings preserve matching performance close to the best observed values while supporting balanced energy transfer, lower mobility cost, higher suitability, and broader Pareto-front exploration. The results show that NSGA-II is relatively stable across the tested hyperparameter ranges, with matching percentages generally remaining within approximately 50.67--53.11\%. Among the crossover settings, a crossover rate of 0.2 achieved the highest average matching performance, with 79.67 matched pairs and a matching rate of 53.11\%, while also producing the highest suitability score of 71.0978. Therefore, the crossover rate was set to 0.2.

For mutation, rates of 0.2 and 0.6 produced the highest average matching performance, both achieving 78.67 matched pairs and a matching rate of 52.45\%. However, mutation rate 0.6 achieved lower mobility cost and higher suitability than mutation rate 0.2 while maintaining the same matching performance. Therefore, the mutation rate was set to 0.6.

For the number of generations, 50 generations achieved the highest matching percentage. However, 250 generations remained within only 0.22 percentage points of this value while achieving the highest transferred energy and the highest suitability score among the tested generation settings. It also produced a substantially lower mobility cost than 50 generations. Since longer runs support broader Pareto-front exploration and improve energy- and suitability-oriented solutions, the number of generations was set to 250.

For population size, 200 achieved the highest matching percentage and suitability score. However, population size 250 remained within 0.45 percentage points of the best matching value while achieving lower mobility cost and higher transferred energy than population size 200. This trade-off was preferred because the final configuration is intended to support balanced multi-objective performance rather than optimize matching percentage alone. Therefore, the population size was set to 250.

Based on this balanced selection rule, the final NSGA-II configuration was set to population size 250, 250 generations, crossover rate 0.2, and mutation rate 0.6.

\begin{table*}[htbp]
\centering

\caption{Sensitivity analysis of NSGA-II hyperparameters for EV--EV energy trading optimization. Results are averaged over three independent tuning runs. Higher values are preferred for matching percentage, transferred energy, and suitability, while lower values are preferred for mobility cost. The dagger ($\dagger$) indicates the selected setting for each hyperparameter.}
\label{tab:sensitivity}
\renewcommand{\arraystretch}{1.15}
\begin{tabular}{llrrrrr}
\toprule
\textbf{Parameter} & \textbf{Value} & \textbf{Avg. $f_1$} & \textbf{Avg. Matching} & \textbf{Avg. $f_2$} & \textbf{Avg. $f_3$} & \textbf{Avg. $f_4$} \\
 &  & \textbf{(matches$\uparrow$)} & \textbf{(\%$\uparrow$)} & \textbf{(energy$\uparrow$)} & \textbf{(mobility$\downarrow$)} & \textbf{(suitability$\uparrow$)} \\
\midrule
Crossover Rate & 0.2$^{\dagger}$ & 79.67 & 53.11 & 1663.39 & 20.0656 & 71.0978 \\
               & 0.4 & 78.00 & 52.00 & 1687.95 & 19.3934 & 69.5478 \\
               & 0.6 & 78.67 & 52.45 & 1674.28 & 19.3675 & 70.6350 \\
               & 0.8 & 78.67 & 52.45 & 1651.21 & 19.6921 & 70.3468 \\
\midrule
Mutation Rate  & 0.2 & 78.67 & 52.45 & 1665.47 & 21.0276 & 69.9304 \\
               & 0.4 & 78.33 & 52.22 & 1661.01 & 20.0535 & 69.9449 \\
               & 0.6$^{\dagger}$ & 78.67 & 52.45 & 1651.21 & 19.6921 & 70.3468 \\
               & 0.8 & 78.33 & 52.22 & 1667.78 & 20.1944 & 69.8514 \\
\midrule
Generations    & 50  & 79.00 & 52.67 & 1660.49 & 21.4623 & 70.0929 \\
               & 100 & 78.67 & 52.45 & 1651.21 & 19.6921 & 70.3468 \\
               & 150 & 78.33 & 52.22 & 1667.26 & 19.2671 & 70.0829 \\
               & 200 & 78.33 & 52.22 & 1670.02 & 18.8385 & 70.1905 \\
               & 250$^{\dagger}$ & 78.67 & 52.45 & 1690.00 & 19.1153 & 70.4040 \\
\midrule
Population Size & 50  & 76.33 & 50.89 & 1626.38 & 20.3394 & 68.3600 \\
                & 100 & 76.00 & 50.67 & 1627.46 & 19.0499 & 68.2846 \\
                & 150 & 78.33 & 52.22 & 1672.19 & 20.7708 & 70.1445 \\
                & 200 & 78.67 & 52.45 & 1651.21 & 19.6921 & 70.3468 \\
                & 250$^{\dagger}$ & 78.00 & 52.00 & 1670.95 & 18.8537 & 69.6477 \\
\bottomrule
\end{tabular}

\vspace{0.5em}
\begin{minipage}{0.96\textwidth}
\footnotesize

\textit{Note:} The selected settings are marked with $\dagger$. They were selected using a balanced multi-objective criterion rather than by optimizing a single objective. The selection prioritizes configurations that maintain competitive matching performance while supporting transferred energy, mobility cost, suitability, and Pareto-front exploration.

\end{minipage}
\end{table*}

\section{Discussion}
\label{sec:discussion}

The experimental results show that the proposed NSGA-II provides a balanced trade-off among the competing objectives considered in EV--EV energy trading. Across both the baseline and state-of-the-art comparisons, NSGA-II consistently achieves higher matching performance, higher transferred energy, and higher suitability than the reference methods, although it does not produce the lowest mobility cost. This pattern is important because it shows that the proposed method is not simply maximizing one metric in isolation. Instead, it identifies compromise solutions that improve the main system-level objectives simultaneously.

The Pareto-front analysis further highlights the multi-objective nature of the problem. The representative Pareto solutions show that improving one objective, such as mobility cost or transferred energy, can lead to reductions in other objectives. In particular, the minimum-mobility solution achieves much lower cost, but with a clear loss in matched pairs, transferred energy, and suitability. This confirms that evaluating EV--EV trading performance using a single objective can be misleading. A balanced solution is therefore more appropriate for practical deployment, where multiple operational goals must be satisfied at the same time.

The hypervolume results provide additional evidence of the effectiveness of NSGA-II. The steady increase and eventual stabilization of hypervolume indicate reliable convergence toward a stable Pareto front, while the larger hypervolume relative to the comparison methods shows that NSGA-II covers a broader and higher-quality trade-off region. This broader coverage is useful in practice because it gives system operators greater flexibility to select solutions according to policy or operational priorities.

The slot-wise analysis also reveals that matching performance varies substantially across time slots. Some slots achieve high matching percentages for both consumers and providers, whereas others remain highly constrained despite the presence of eligible participants. This variation indicates that temporal feasibility and counterpart availability strongly influence final outcomes. In this setting, the benefit of NSGA-II lies in its ability to exploit favorable slots while preserving feasibility under more constrained temporal conditions.

In summary, the results suggest that the main advantage of NSGA-II is its ability to balance competing objectives in a dynamic and constrained environment. Rather than minimizing mobility cost alone or relying on locally attractive assignments, it produces solutions that better align with the broader goals of EV--EV energy trading, namely effective matching, high energy transfer, and strong assignment suitability.

\section{Summary of Findings}
\label{sec:summary}

The results show that NSGA-II achieves the balanced performance across matching coverage, transferred energy, and charging-node suitability, although this comes with higher mobility cost than some comparison methods. This mobility cost should therefore be interpreted as an explicit trade-off rather than as an independently optimized outcome. The Pareto-front and hypervolume analyses further confirm that the proposed method provides broader and higher-quality trade-off solutions than the comparison methods, allowing different solutions to prioritize matching coverage, transferred energy, mobility cost, or charging-node suitability according to operational needs. In addition, the slot-wise analysis shows that performance depends strongly on temporal opportunity and feasibility, which further supports the need for coordinated multi-objective optimization in EV--EV energy trading systems.

\section{Conclusion and Future Work}
\label{sec:conclusion}

This study presented a prediction-guided, multi-objective optimization framework for provider--consumer matching in EV--EV energy trading. The framework integrates prediction outputs from EVNextTrade with a multi-objective MILP formulation and an NSGA-II search procedure, jointly considering matching coverage, transferred energy, mobility cost, and charging-node suitability under spatial, temporal, and charging node capacity constraints. By leveraging prediction-derived suitability scores as guidance signals, the framework improves coordination among decentralized EV participants without relaxing the underlying operational constraints.

Experimental results confirmed the effectiveness of the proposed approach against both heuristic and learning-based baselines and against state-of-the-art EV--EV matching methods. NSGA-II achieved the balanced performance, simultaneously improving matching coverage, transferred energy, and suitability relative to DAMA and BMNN, while accepting higher mobility cost as an explicit trade-off. Pareto-front and hypervolume analyses further demonstrated reliable convergence behaviour and a wide range of feasible trade-offs, supporting flexible solution selection across different operational priorities.

Beyond the empirical improvements, this study contributes a practical decision-support framework for coordinated EV--EV energy trading under spatiotemporal uncertainty. By combining prediction-guided candidate evaluation with multi-objective optimization, the proposed approach moves beyond single-criterion matching and enables more informed, efficient, and flexible coordination among mobile EVs. This is particularly relevant for future charging ecosystems, where decentralized energy exchange must operate under congestion, heterogeneous user conditions, and dynamically changing infrastructure availability.


Several limitations remain. The current formulation treats recommendation outputs as fixed inputs and assumes stationary charging nodes. Although the formulation enforces event-level energy feasibility through generated demand, supply, journey-feasibility, and transferable-energy constraints, repeated decision events generated by the same physical EV are treated as independent trading opportunities rather than as a longitudinal sequence of post-trade battery-state updates. In addition, the evaluation is conducted using a single-day urban dataset, which limits the ability to capture broader temporal variation and dynamic operating conditions. Building on these observations, future work will proceed in four directions. First, tighter integration between prediction and optimization will be investigated, including online or closed-loop frameworks in which predictions and matching decisions are updated jointly over time. Second, as EV participation increases, scalability and privacy considerations will be explored through distributed, federated, or privacy-preserving optimization strategies. Third, the framework will be extended by incorporating physical-EV-level battery-state propagation across repeated decision events. Finally, the framework will be extended to more dynamic charging settings, such as mobile charging units, charging-as-a-service platforms, and direct EV--EV energy exchange without fixed intermediaries.

\section*{Acknowledgments}
The authors acknowledge that this research is partly supported by the QUT Postgraduate Research Award (QUTPRA)
scholarship.


\bibliographystyle{IEEEtran}
\bibliography{references}













\section{Biography Section}
 


\begin{IEEEbiography}[{\includegraphics[width=1in,height=1.25in, clip,keepaspectratio]{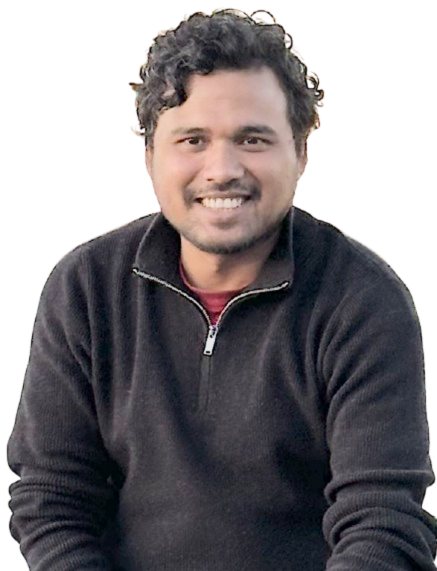}}]{Md Mahfujur Rahman} (Graduate Student Member, IEEE) has received both a B.Sc. (Hons.) and an M.Sc. in Information Technology from Jahangirnagar University, Savar, Dhaka, Bangladesh. With extensive practical experience and research expertise, he has successfully completed over 15 IT projects and authored 26 research publications, showcasing his dedication to advancing IT knowledge and research. Rahman was also involved in the development and management of the conference software PROCONF, which is now recognised worldwide. Rahman is a Cisco Certified Network Associate (CCNA) and holds certifications in Database Design and Implementation with Web Applications. He began his career as a software engineer before transitioning to academia as a faculty member. His leadership skills were evident in his role as the former General Secretary of the IEEE Student Branch at JU. Rahman's research encompasses machine learning, deep learning, optimization, and digital ledger technology. Currently, he is pursuing a PhD in EV--EV energy trading at the School of Information Systems, Queensland University of Technology, aiming to make significant contributions to sustainable energy solutions. His journey reflects an unwavering dedication to IT and a commitment to advancing the boundaries of knowledge and technology.
\end{IEEEbiography}
\vspace{-1mm}
\begin{IEEEbiography}[{\includegraphics[width=1in,height=1.25in, clip,keepaspectratio]{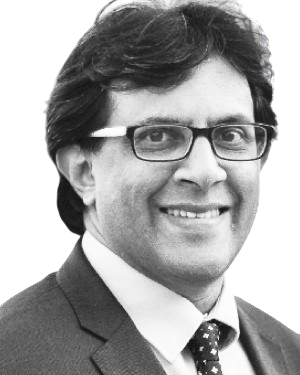}}]{Alistair Barros} is Head of School and Academic Program Leader of Service Science at QUT’s School of Information Systems at the Faculty of Science, QUT. He has a PhD from the University of Queensland and ICT experience across academic and industry organisations, including being Global Research Leader at SAP AG. Alistair's research focus is on the design, evolution, interoperability, and optimization of enterprise systems through contemporary cyber-physical settings, enabled by Cloud, Industrial Internet of Things, and Blockchain platforms. His research interests include service computing methods and techniques applied to: software architectures and microservices; business process management; model-based systems re-engineering; and distributed service optimization and coordination. His work previously applied to service industries in banking and the public sector, and has extended to cyber-physical domains including construction, manufacturing, and supply chains. Alistair has published 180+ articles, which include six edited books, 130 peer-reviewed journals, conference, and book chapter articles. He also has 17 filed US patents. In terms of research publication metrics, he has a h-index of 40, and his articles have attracted 11,000 citations, according to Google Scholar. Among his publication highlights are: the “Workflow Patterns” article, which is the most cited in the Business Process Management (BPM) field; and the “Service Interaction Patterns” which won the “Test of Time” award at international BPM 2015 conference for highest impact BPM 2005-6 paper over ten years. Alistair has led the  Australian Research Council, the Cooperative Research Centre, and the EU Framework.
\end{IEEEbiography}
\begin{IEEEbiography}[{\includegraphics[width=1in,height=1.25in, clip,keepaspectratio]{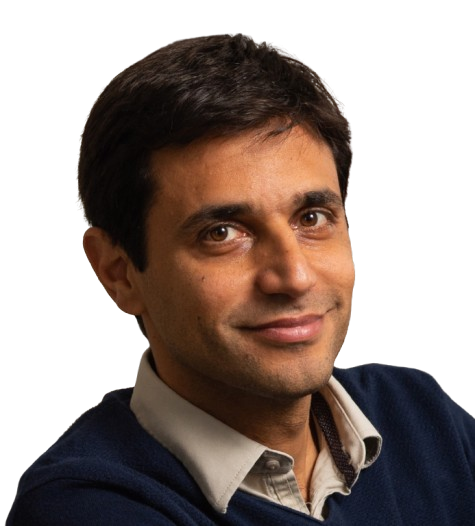}}]{Raja Jurdak} (Senior Member, IEEE)
received the B.E. degree in computer and communications engineering from the American University of Beirut, Beirut, Lebanon, in 2000, the M.S. degree in computer networks and distributed computing from the Department of Electrical and Computer Engineering, University of California at Irvine, Irvine, CA, USA, in 2001, and the Ph.D. degree in information and computer science from the University of California at Irvine, Irvine, in 2005., He is a Professor of Distributed Systems and the Chair of Applied Data Sciences with Queensland University of Technology, Brisbane, QLD, Australia, where he is also the Director of the Trusted Networks Laboratory. He is also an Adjunct Professor with the University of New South Wales, Sydney, NSW, Australia, and a Visiting Researcher with Data61, CSIRO, Brisbane. Previously, he established and led the Distributed Sensing Systems Group, Data61, CSIRO. He also spent time as a Visiting Academic with the Massachusetts Institute of Technology, Cambridge, MA, USA, and Oxford University, Oxford, U.K., in 2011 and 2017, respectively. He has published over 250 peer-reviewed publications, including three authored books on IoT, blockchain, and cyber-physical systems. His publications have attracted over 14700 citations, with an H-index of 54. His research interests include trust, mobility, and energy efficiency in networks.  Prof. Jurdak serves on the organizing and technical program committees of top international conferences, including Percom, ICBC, IPSN, WoWMoM, and ICDCS. He serves on the editorial boards of Ad Hoc Networks and Scientific Reports (Nature). He is an IEEE Computer Society Distinguished Visitor.
\end{IEEEbiography}
\begin{IEEEbiography}[{\includegraphics[width=1in,height=1.25in, clip,keepaspectratio]{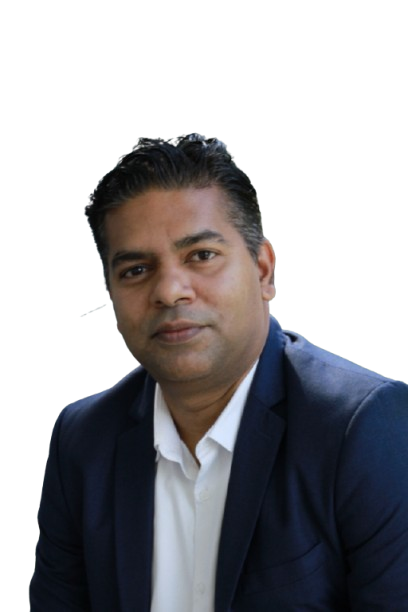}}]{Darshika Koggalahewa} 
Dr. Darshika Koggalahewa is a leading expert in Enterprise Systems, Data Science, and Artificial Intelligence. He holds a PhD in Data Science and Machine Learning from QUT, where he received the Faculty of Science Best Thesis Award. As Subject Area Lead in Enterprise Systems at QUT, Dr. Koggalahewa leads curriculum innovation and applied research in analytics, software engineering, and enterprise computing, partnering with SAP and Microsoft to deliver industry-aligned educational programs. His teaching and academic leadership have been recognised through multiple awards, including QUT-IS’s Best Individual Teaching Award, Educator of the Year (Runner-up, 2023), Best Early Career Award, and Outstanding Achievements Award. Dr. Koggalahewa’s research interests include Enterprise Data Management, Data Science, Machine Learning, and Industrial IoT (IIoT), focusing on the optimization of large-scale business processes. His expertise spans Enterprise Data Management and Machine Learning, aligning data governance, architecture, and analytics to enable AI/ML-driven enterprise systems. He collaborates with CSIRO and the Australian logistics industry to optimise port logistics through data-driven enterprise solutions—delivering real-world impact through innovation. Dr. Koggalahewa has published extensively in Data Science, Machine Learning, and Enterprise Computing and has led the development of undergraduate and postgraduate programs in Australia and Sri Lanka. He has been instrumental in capacity building within the Sri Lankan education sector, designing and delivering IT and Software Engineering programs, while also working closely with the Sri Lankan software industry to advance skills in enterprise solution development and mobile computing.
\end{IEEEbiography}
%




\vfill

\end{document}